\documentclass[12pt,french]{article}
\usepackage{geometry}
\usepackage{amsmath}\usepackage{amssymb}
\usepackage{amsthm}
\usepackage{graphicx}
\usepackage{stmaryrd}
\usepackage[french, english]{babel}
\usepackage[T1]{fontenc} 
\usepackage[babel=true]{csquotes}
\usepackage[backend=biber,
style=alphabetic,
giveninits=true, maxbibnames=99, maxalphanames=99, url=false, doi=false, isbn=false, date=year]{biblatex}
\DeclareFieldFormat{postnote}{#1} 
\AtEveryBibitem{\clearfield{issue}\clearfield{url}\clearfield{urlyear}
    \clearfield{urlmonth}\clearfield{pubstate}}
\DeclareFieldFormat{pages}{#1}
\renewbibmacro{in:}{}
\usepackage{xifthen}
\usepackage{tikz-cd}
\usepackage{multirow}
\usepackage{caption}

\newcommand{\nc}{\newcommand}
\nc{\nt}{\newtheorem}
\nc{\frg}{\mathfrak{g}}
\nc{\frsl}{\mathfrak{sl}}
\nc{\rhohat}{\hat{\rho}}
\nc{\id}{\mathrm{id}}
\nc{\Thmdot}[2]{{#1}.~{#2}}

\nc{\fru}{\mathfrak{u}}
\nc{\du}{\overdot{\fru}}
\nc{\be}{\bold{e}}
\nc{\bes}{\be^{\times}}
\nc{\cG}{\check{G}}
\nc{\Ed}{\mathrm{E}_2}
\nc{\CC}{\mathbb{C}}
\nc{\II}{\mathbb{I}}
\nc{\hvee}{h^{\vee}}

\nc{\LG}{{^LG}}
\nc{\Iuminus}{I^{u,-}}
\nc{\tcN}{\tilde{\cN}}
\nc{\Lp}{\mathrm L}
\nc{\dR}{_{\mathrm{dR}}}
\nc{\Gsc}{G^{\mathrm{sc}}}
\nc{\red}{_{\mathrm{red}}}
\nc{\disj}{\mathrm{disj}}
\nc{\nilp}{\mathrm{nilp}}
\nc{\Ran}{\mathrm{Ran}}
\nc{\Rannv}{\mathrm{Ran}^{\mathrm{nv}}}
\nc{\DRan}[1][]{\ifthenelse{\isempty{#1}}{D_{\Ran}}{D_{\Ran_{#1}}}}
\nc{\cCRan}{\cC_{\Ran}}
\nc{\DcRan}[1][]{\ifthenelse{\isempty{#1}}{\Dc_{\Ran}}{\Dc_{\Ran_{#1}}}}
\nc{\Rani}{\mathrm{Ran}^{\subseteq}}
\nc{\Sing}{\mathrm{Sing}}
\nc{\Arc}{\mathrm{L}^+}
\nc{\LpG}{\Lp G}
\nc{\AG}{\Arc G}
\nc{\Lh}{h^{\vee}}
\nc{\Lg}{{^L\mathfrak{g}}}
\nc{\SL}{\mathrm{SL}}
\nc{\Spin}{\mathrm{Spin}}
\nc{\GL}{\mathrm{GL}}
\nc{\PGL}{\mathrm{PGL}}
\nc{\rmA}{\mathrm{A}}
\nc{\rmB}{\mathrm{B}}
\nc{\rmC}{\mathrm{C}}
\nc{\rmT}{\mathrm{T}}
\nc{\rmE}{\mathrm{E}}
\nc{\rmF}{\mathrm{F}}
\nc{\rmN}{\mathrm{N}}
\nc{\rmG}{\mathrm{G}}
\nc{\sph}{\mathrm{sph}}
\nc{\Psph}{P^{\sph}}

\nc{\rmD}{\mathrm{D}}
\nc{\Dc}{D^{\times}}
\nc{\Gr}{\mathrm{Gr}}
\nc{\Grcirc}{\mathrm{Gr}_G^{\circ}}
\nc{\GrRan}[1][G]{\mathrm{Gr}_{#1,\Ran}}
\nc{\GrRanJ}[2][]{
\ifthenelse
{\isempty{#1}}{\Gr_{#2,\Ran}}{\Gr_{#2,\Ran_{#1}}}
}
\nc{\Grlambda}{\Gr_{\Gsc,\Ran_J}^{\boxtimes\lambda_j}}
\nc{\GrRancirc}{\mathrm{Gr}_{G,\Ran}^{\circ}}

\nc{\sunderl}[1]{s_{\underline{#1}}}
\nc{\Normalizer}[2]{N_{#1}(#2)}

\nc{\AAA}{\mathbb{A}}
\nc{\OO}{\mathbb{O}}
\nc{\GG}{\mathbb{G}}
\nc{\PP}{\mathbb{P}}
\nc{\ZZ}{\mathbb{Z}}

\nc{\cC}{\mathcal{C}}
\nc{\cA}{\mathcal{A}}
\nc{\cE}{\mathcal{E}}
\nc{\cF}{\mathcal{F}}
\nc{\cL}{\mathcal{L}}
\nc{\cN}{\mathcal{N}}
\nc{\cO}{\mathcal{O}}
\nc{\cV}{\mathcal{V}}

\nc{\ttt}{\text{-}}
\nc{\Gm}{\GG_m}
\nc{\Ga}{\GG_a}
\nc{\frga}{\frg_a}
\nc{\Fib}{\mathrm{Fib}}
\nc{\Indcoh}{\mathrm{Indcoh}}
\nc{\Qcoh}{\mathrm{Qcoh}}
\nc{\Coh}{\mathrm{Coh}}
\nc{\Dmod}{\rmD\ttt\mathrm{mod}}
\nc{\Ind}{\mathrm{Ind}}
\nc{\Perf}{\mathrm{Perf}}
\nc{\Sym}{\mathrm{Sym}}
\nc{\Rep}{\mathrm{Rep}}
\nc{\act}{\mathrm{act}}
\nc{\Mod}{\mathrm{Mod}}
\nc{\tMod}{\ttt\Mod}
\nc{\Hom}{\mathrm{Hom}}
\nc{\IC}{\mathrm{IC}}
\nc{\ad}{\mathrm{ad}}
\nc{\Ad}{\mathrm{Ad}}
\nc{\rmLS}{\mathrm{LS}}
\nc{\Spf}{\mathrm{Spf}}
\nc{\Spec}{\mathrm{Spec}}
\nc{\cExp}{\mathcal{Exp}}
\nc{\pr}{\mathrm{pr}}
\nc{\Al}{\mathrm{A}_{\ell}}
\nc{\Bl}{\mathrm{B}_{\ell}}
\nc{\Cl}{\mathrm{C}_{\ell}}
\nc{\Dl}{\mathrm{D}_{\ell}}
\nc{\El}{\mathrm{E}_{\ell}}

\nc{\Vect}{\mathrm{Vect}}
\nc{\tonte}{^{\fatslash}}

\nc{\ito}{\overset{\sim}{\to}}

\nc{\frus}{\mathfrak{u}_{\sigma}}
\nc{\Us}{U_{\sigma}}
\nc{\As}{\cA_{\sigma}}
\nc{\Ls}{\cL_{\sigma}}
\nc{\Fs}{\cF_{\sigma}}
\nc{\Vs}{\cV_{\sigma}}

\nc{\LsZRan}[1][]{
\ifthenelse
{\isempty{#1}}{\cL'_{\sigma,\Ran}}{\cL'_{\sigma,\Ran_{#1}}}
}
\nc{\Lsk}[1][\chi]{\cL_{\sigma,#1}}
\nc{\Fsk}[1][\chi]{\cF_{\sigma,#1}}
\nc{\Vsk}[1][\chi]{\cV_{\sigma,#1}}
\nc{\LskRan}[2][]{
\ifthenelse
{\isempty{#1}}{\cL_{\sigma,#2,\Ran}}{\cL_{\sigma,#2,\Ran_{#1}}}
}
\nc{\FskRan}[2][]{
\ifthenelse
{\isempty{#1}}{\cF_{\sigma,#2,\Ran}}{\cF_{\sigma,#2,\Ran_{#1}}}
}
\nc{\VskRan}[2][]{
\ifthenelse
{\isempty{#1}}{\cV_{\sigma,#2,\Ran}}{\cV_{\sigma,#2,\Ran_{#1}}}
}
\nc{\FskRannv}[2][]{
\ifthenelse
{\isempty{#1}}{\cF_{\sigma,#2,\Rannv}}{\cF_{\sigma,#2,\Rannv_{#1}}}
}
\nc{\FsRan}{\cF_{\sigma,\Ran}}
\nc{\VsRan}{\cV_{\sigma,\Ran}}
\nc{\Ask}{\cA_{\sigma,\chi}}
\nc{\Sls}{\mathrm{Sl}_{\sigma}}
\nc{\SOk}{S_{O,\chi}}
\nc{\Ms}{M_{\sigma}}
\nc{\Msn}{M_{\sigma}^{\nilp}}
\nc{\OMs}{O_{\Ms}}
\nc{\OMst}{O_{\Ms}\tonte}
\nc{\omegat}{\omega\tonte}
\nc{\Zs}{Z_{\sigma}}
\nc{\Zt}{Z_{\tau}}
\nc{\bphi}{\bar{\phi}}
\nc{\Sc}{S_{\mathrm coh}}
\nc{\Aphi}{\AAA^1_{\phi}}
\nc{\Aphid}{\AAA^{1,*}_{\phi}}
\nc{\Gaphi}{\GG_{a,\phi}}
\nc{\frgaphi}{\frg_{a,\phi}}
\nc{\frgaphid}{\frg^*_{a,\phi}}
\nc{\Lk}[1][\kappa]
{L_{#1}}
\nc{\LkL}[2][\kappa]{L_{#1}[#2]}
\nc{\wlambda}{w_{\lambda}}
\nc{\wthree}[1]{w(#1)}

\nc{\Hk}{\mathrm{Hk}}
\nc{\Sph}{\mathrm{Sph}}
\nc{\SphG}{\Sph_G}
\nc{\HkIG}{\Hk^I_G}
\nc{\FibxU}{\Fib_G^{x,U}}
\nc{\Fibx}{\Fib_G^{\hat{x}}}
\nc{\Fibhxy}{\Fib_G^{\hat{x};y,B^-}}
\nc{\FibxBy}{\Fib_G^{x,B;y,B^-}}
\nc{\FibxUy}{\Fib_G^{x,U;y,B^-}}
\nc{\HkGmonodrom}{\Dmod\bigl(I
\rotatebox[origin=c]
{30}{\textbrokenbar} G(\Dc)
\rotatebox[origin=c]
{-30}{\textbrokenbar}I\bigr)}
\nc{\Hks}{\mathrm{Hk}^{\mathrm{spec}}}
\nc{\Hkst}[1][\Ga]{\mathrm{Hk}^{\mathrm{spec},\fatslash}_{#1}}
\nc{\HkstRan}[1][\Ga]{\mathrm{Hk}^{\mathrm{spec},\fatslash}_{#1,\Ran}}
\nc{\Sphs}{\mathrm{Sph}^{\mathrm{spec}}}
\nc{\Sphst}{\mathrm{Sph}^{\mathrm{spec},\fatslash}}
\nc{\Sat}{\mathrm{Sat}^{\mathrm{spec}}}
\nc{\SphsLG}[1][\LG]{\Sphs_{#1}}
\nc{\SphsLGRan}[1][]{
\ifthenelse
{\isempty{#1}}{\Sphs_{\LG,\Ran}}{\Sphs_{\LG,\Ran_{#1}}}
}
\nc{\SphsZRanJ}[2][]{
\ifthenelse
{\isempty{#1}}{\Sphs_{#2,\Ran}}{\Sphs_{#2,\Ran_{#1}}}
}
\nc{\SphstZRan}[1]{\Sphst_{#1,\Ran}}
\nc{\SphsZRan}[1]{\Sphs_{#1,\Ran}}
\nc{\SphGRan}[1][]{
\ifthenelse
{\isempty{#1}}{\Sph_{G,\Ran}}{\Sph_{G,\Ran_{#1}}}
}
\nc{\HksGa}{\Hks_{\Ga}}
\nc{\SatLG}{\Sat_{\LG}}
\nc{\cERan}{\cE_{\Ran}}
\nc{\LS}[2][\LG]{\rmLS_{#1}^{#2}}
\nc{\LSUs}{\rmLS_{\Us\subset\LG\times\Ga}^{\Dc\subset D}}
\nc{\LSZUs}{\rmLS_{\Zs\Us\subset\Zs\times\LG\times\Ga}^{\Dc\subset D}}
\nc{\LSZUsRan}{\rmLS_{\Zs\Us\subset\Zs\times\LG\times\Ga}^{\DcRan\subset\DRan}}
\nc{\tlambda}{t^{\lambda}}
\nc{\tlambdaj}{t^{\boxtimes\lambda_j}}
\nc{\Lhat}{\hat{\Lambda}}
\nc{\Deltaw}[1][w]{\Delta_{#1}}
\nc{\Deltavi}{\Delta_{v^{-1}}}
\nc{\nablavi}{\nabla_{v^{-1}}}
\nc{\Cvi}{C_{v^{-1}}}
\nc{\Deltatilde}{\tilde{\Delta}}
\nc{\nablatilde}{\tilde{\nabla}}
\nc{\Deltawi}[1][w]{\Delta_{{#1}^{-1}}}
\nc{\ICOk}{\mathrm{IC}_{O,\chi}}
\nc{\kp}{\kappa^+}

\newtheorem*{thm*}{Théorème}
\newtheorem*{conj*}{Conjecture}
\newtheorem*{question*}{Question}

\newtheorem{cond}{Condition}

\theoremstyle{remark}

\title{L'Algèbre de Hecke chirale nilpotente}
 \author{Lizao Ye}
\date{}
\begin{document}

\maketitle
\selectlanguage{french}
\begin{abstract}
    Nous étendons aux orbites nilpotentes la notion d'algèbre de Hecke chirale introduite par Beilinson-Drinfeld. En analysant leurs composantes isotypiques, nous produisons de nombreux modules inédits sur des algèbres vertex affines simples à des niveaux entiers non-admissibles, ainsi que leurs façons de fusionner.
\end{abstract}
\selectlanguage{english} 
\begin{abstract}
We extend to nilpotent orbits the notion of chiral Hecke algebra introduced by Beilinson-Drinfeld. Upon analysing their isotypic components, we produce many new modules over simple affine vertex algebras at non-admissible integer levels, as well as their ways of fusion.
\end{abstract}

\selectlanguage{french} 
Récemment, des algèbres vertex affines liées à certaines orbites nilpotentes apparaissent dans des contextes 
variés, de façon assez sporadique. Il n'est pas clair au premier abord comment les généraliser aux autres orbites nilpotentes. Cependant, dans un contexte différent, celui de la correspondance de Langlands géométrique, des objets naturels liés à toutes les orbites nilpotentes, quoique du groupe dual de Langlands, ont bien été construits. Ce sont les \textit{faisceaux de Whittaker spectraux généralisés}. Le but de cet article est de montrer un lien entre ces deux. L'idée est que l'application de Spaltenstein permettrait de revenir aux orbites nilpotentes du groupe original.

En effet, les algèbres vertex affines sont les objets factorisables parmi les représentations d'algèbres de Lie affines. Une façon naturelle de les produire est de prendre les sections globales --- tordues par les fibrés en droites factorisables --- des D-modules factorisables sur la grassmannienne affine factorisable de Beilinson-Drinfeld, qui ont construit leur algèbre de Hecke chirale exactement par cette méthode. Nos faisceaux spectraux, étant factorisables, donnent de tels D-modules via l'équivalence de Satake dérivée factorisable. En fait, leur construction peut être vue comme le cas correspondant à l'orbite zéro parmi les nôtres, d'où le titre de cet article.

Il se trouve que ces faisceaux spectraux admettent l'action supplémentaire du quotient réductif du centralisateur d'un élément de l'orbite nilpotente et se décomposent de la même manière que la représentation régulière de ce groupe. Nous montrons que cette action est compatible avec la structure de factorisation, ce qui implique une décomposition de nos algèbres vertex dont les composantes se fusionnent tout comme les représentations irréductibles de ce groupe. Un fait important est que chaque D-module isotypique sous-jacent est supporté sur une seule composante connexe de la grassmannienne affine, si bien que l'on peut les transporter tous sur la composante neutre de manière cohérente et n'utiliser que les fibrés en droites factorisables là-dessus pour la section globale tordue. C'est moins restreignant que d'employer les fibrés factorisables sur la grassmannienne affine tout entière puisque là il faut compter en plus le mélange des différentes composantes connexes. 
Finalement on va comparer les algèbres vertex ainsi construites avec la littérature récente, suivant laquelle on supposera que nos fibrés en droites
soient de niveau plus grand que critique. C'est ici qu'on a recours à la correspondance de Langlands géométrique transformant les faisceaux spectraux aux D-modules automorphes, les sections globales tordues desquels produisant les mêmes algèbres vertex grâce au fait que ces D-modules sont liés aux ceux sur la grassmannienne affine dont on vient de parler par la transformée de Radon longue. En fait, on va donner une description de ces D-modules automorphes en termes des faisceaux pervers simples. Lorsque l'orbite nilpotente est la régulière, c'est-à-dire la plus grande, le D-module automorphe est tout simplement le faisceau constant. Si l'orbite est proche de la régulière, les faisceaux pervers simples intervenants seront \enquote{proches} du faisceau constant. Dans ces cas, leurs sections globales tordues à un niveau convenable plus grand que critique sont contrôlées, ce qui aboutirait à des formules explicites de nos algèbres vertex permettant la comparaison avec la littérature.

Les niveaux qui apparaissent dans ces articles récents sont d'ailleurs plus petits que le niveau trivial. Ils sont donc \enquote{singuliers}, et non-admissible dans la terminologie de Kac-Wakimoto. Cette restriction n'est pas nécessaire dans notre construction générale. Toutefois, on en a besoin, si l'on veut une formule explicite, lorsque l'orbite nilpotente n'est pas la régulière. C'est parce que notre description des objets automorphes contient une ambiguïté dans ces cas, qui ne pourrait être tuée que par la section globale tordue par un poids singulier. Lorsque le poids est de plus $\rhohat$-dominant, on sait que la partie tuée est engendrée par des faisceaux pervers simples déterminés par la singularité du poids. Nous allons exécuter des transformées de Radon \enquote{courtes} pour nous ramener à cette situation favorable si l'on n'en est pas loin au départ, autrement dit si le niveau est proche de zéro. La singularité du niveau, un défaut, devient ainsi un atout. Nous produisons de cette manière une large famille des modules irréductibles sur les algèbres vertex affines simples à un niveau entier non-admissible ainsi que leurs lois de fusion. Ces résultats sont pour la plupart inédits. Nous obtenons par exemple:
\begin{align*}
    \Lk[-2](\SL_3)\tMod &\supseteq\Rep(\SL_3) ,\\\Lk[-3](\rmG_2)\tMod &\supseteq\Rep(\SL_2). 
\end{align*}
Des tableaux plus détaillés se trouvent vers la fin. On y voit aussi quelques cas déjà connus tels que \begin{align*}
    \Lk[-1](\SL_n)\tMod &\supseteq\Rep(\Gm), n\geqslant 3,\\\Lk[-2](\rmG_2)\tMod &\supseteq\Rep(S_3) 
\end{align*}  construits plus tôt dans~\cite{adamovicFusionRulesComplete2014} et~\cite{adamovicSemisimplicityModuleCategories2024} respectivement, à travers des méthodes complètement différentes. Ils concluent de plus que ces deux inclusions sont en fait des équivalences.
Lorsque l'orbite nilpotente est la régulière, le D-module automorphe correspondant est le faisceau constant et exige un niveau régulier, autrement dit, dominant. L'algèbre vertex associée est alors une somme des courants simples indexés par les composantes connexes, lesquels étaient discutés pour la première fois dans~\cite{fuchsConnectionWZWFree1987}.

Le lecteur aurait dû se rendre compte que notre méthode ne donne des formules explicites que pour les orbites proches de la régulière et pour les niveaux plus grands que critique. On doit cependant mentionner le cas de l'orbite zéro. Ici le D-module associé sur la grassmannienne affine est déjà une somme explicite des faisceaux pervers simples. Sa section globale tordue à un niveau plus petit que critique est alors explicitement une somme des représentations irréductibles de l'algèbre de Lie affine et donne justement l'algèbre de Hecke chirale de Beilinson-Drinfeld sans passer au côté automorphe. Pour un compte rendu, voir~\cite{shapiroBRSTReductionChiral2009}.

\textbf{Conventions.} Le corps de base est le corps des nombres complexes. Les catégories sont dérivées et ind-complétées. Les foncteurs sont dérivés. L'image réciproque signifie la $!$-version. Les catégories de D-modules sont renormalisées. Cela correspond à ne pas imposer des restrictions sur le support singulier au côté cohérent. Les D-modules monodromiques sont de monodromie unipotente. La numérotation des sommets du diagramme de Dynkin d'un groupe presque-simple suit~\cite[Table Fin, p. 53]{kacInfiniteDimensionalLieAlgebras1990}. Le poids fondamental pour la sommet~$i$ est noté~$\Lambda_i$. On y ajoute la sommet~$0$ pour former le diagramme de Dynkin affine. Les poids fondamentaux affines sont notés~$\Lhat_i$. Le groupe de Weyl affine, sans l'épithète \enquote{étendu}, signifie la version non-étendue. Les mots \enquote{positif}, \enquote{négatif} n'excluent pas zéro.

\textbf{Remerciements.} Je suis reconnaissant à de nombreuses personnes, sans qui ce projet n'aurait pas pu voir le jour. Vincent Lafforgue l'a appuyé depuis fin~2018, quand je n'en avais qu'une vision provisoire. Il avait une idée similaire. Vadim Schechtman s'y est très intéressé 
et a organisé une rencontre à Toulouse l'année suivante. Michael Finkelberg m'a appris les faits sur la localisation. Toshiyuki Tanisaki m'a mentionné l'article~\cite{kashiwaraKazhdanLusztigConjectureAffine1996}. Sam Raskin m'a mentionné le traitement du faisceau automorphe constant dans~\cite{ben-zviRelativeLanglandsDuality2024}. Avant tout cela, Sergey Lysenko m'a mentionné début 2016 la trialité du faisceau automorphe minimal pour~$\rmD_4$ ainsi que son lien avec~$\rmG_2$.

Ce projet était conçu à l'Institut Élie Cartan de Lorraine et a pris forme pendant les séjours successifs de l'auteur à l'Université nationale de recherche \enquote{École des hautes études en sciences économiques}, Moscou, au Centre Morningside des Mathématiques, l'Académie chinoise des sciences, Pékin, et à l'Université Aalto, Helsinki. Je remercie tous ces instituts pour leurs supports.

Les faisceaux spectraux sont introduits dans l'article récent~\cite{ben-zviRelativeLanglandsDuality2024} de Ben-Zvi, Sakellaridis et Venkatesh. Notre exposition dans cette partie suit largement leur article.  Soit~$G$ un groupe réductif, $\LG$ son groupe dual au sens de Langlands. Notons $\frg,\Lg$ leurs algèbres de Lie respectives. Les faisceaux de Whittaker spectraux généralisés font partie de la famille plus grande des L-faisceaux, tous vivent dans le côté spectral de la correspondance de Langlands qu'on identifie pour le moment avec la  catégorie de Satake dérivée spectrale $\SatLG=\Ind\Bigl(\Perf\left(\Sym(\Lg [-2])/\LG\right)\Bigr)$ via la dualité de Koszul. Plus loin on aura à revenir à sa version originale en termes de systèmes locaux pour voir, en fait déjà pour formuler, la propriété factorisable de ces faisceaux. 

En général, les L-faisceaux sont associés aux données suivantes. Il faut d'abord une variété symplectique~$M$ munie d'une action de $\LG\times\Gm$ telle que la forme symplectique est préservée par $\LG$ et varie de poids~$-2$ sous la poussée vers l'avant de $\Gm$. Autrement dit, le champ de bivecteurs associé est invariant sous  $\LG$ et varie de poids~2 sous $\Gm$. Cette reformulation permet d'inclure aussi les variétés de Poisson. Le choix du poids dicte en gros que le crochet de Poisson augmente de~2 la graduation sur les fonctions. On voit aisément que la variété de Poisson $\Lg^*$ munie de l'action scalaire au carré de $\Gm$ vérifie cette condition. Cet exemple sert à exprimer l'autre part des données. C'est une application moment $\mu: M\to\Lg^*$ qui est $\LG\times\Gm$-équivariante. 
 
S'il ne pose pas de confusion, on parlera simplement du L-faisceau associé à~$M$ pour désigner celui associé aux données ci-dessus. Il est presque $\mu_* O_M$ où $O_M$ est le faisceau structural sur~$M$. Le problème est que ce dernier se trouve dans la catégorie $\Qcoh(\Lg^*/\LG)$, ou, au mieux, est un objet $\Gm$-équivariant dedans.  Il existe néanmoins une procédure canonique, figurément appelée \enquote{tonte}, pour passer à la vraie catégorie de Satake dérivée spectrale. Laissons-nous la rappeler brièvement puisqu'elle joue un rôle important.
 
La catégorie des faisceaux quasi-cohérents sur $\Gm$ est monoïdale sous la convolution. Supposons donnée une catégorie~$\cC$ sous son action. Les objets équivariants, c'est-à-dire les objets~$c$ munis des isomorphismes compatibles  $i_V: V\star c \ito V\otimes c$ pour tout $V\in\Qcoh(\Gm)$ où~$\star$ désigne l'action et $\otimes$ désigne la tensorisation avec l'espace vectoriel banal sous-jacent, forment une catégorie, notée $\cC^{\Gm}$,  sur laquelle $\Rep(\Gm)$ agit en ne modifiant que les isomorphismes: un caractère $\lambda$, pouvant être considéré comme élément de $O_{\Gm}$, envoie $i_V$ à $(\lambda\otimes 1)\circ i_V$. On peut retrouver la catégorie initiale $\cC$ en trivialisant cette action. La tonte consiste à  munir $\cC^{\Gm}$ d'une nouvelle action de $\Rep(\Gm)$: l'effet $\act_{\lambda}$ de $\lambda$ est remplacée par $\act_{\lambda}\circ [-n]$ lorsque $\lambda$ désigne la $n$-ième puissance. On trivialise cette action décalée pour arriver à la version tondue de~$\cC$, notée $\cC\tonte$. Un objet $c\in\cC$, si équivariant, peut être considéré comme objet de $\cC^{\Gm}$; son image dans $\cC\tonte$, notée $c\tonte$, est la version tondue de~$c$. 
Remarquons que pour arriver à $\cC\tonte$, on pourrait aussi procéder autrement en changeant non pas l'action de $\Rep(\Gm)$  mais la façon de la trivialiser. En effet, on peut prendre la nouvelle trivialisation  $\Rep(\Gm)\to\Vect$ qui est la composition de l'automorphisme monoïdal de  $\Rep(\Gm)$ envoyant~$V$ à  $V\tonte :=\oplus V^i[i]$ où $V^i$ désigne la composante de poids~$i$ de~$V$, et la trivialisation habituelle oubliant la graduation, pour réaliser la même tonte.
   
L'idée est que pour des dg-algèbres graduées $A$ convenables, on devrait avoir $(A\tMod)\tonte\ito A\tonte\tMod$. Par exemple, pour la $\Gm$-variété $\Lg^*$ comme ci-dessus, la version tondue de $\Qcoh(\Lg^*)$ est $\Ind\Bigl(\Perf\left(\Sym(\Lg [-2])\right)\Bigr)$; celle de $\Qcoh(\Lg^*/\LG)$ est $\SatLG$. On peut maintenant définir le L-faisceau associé à~$M$ comme la version tondue de  $\mu_* O_M$.
 
Les faisceaux de Whittaker spectraux généralisés sont des L-faisceaux associés aux orbites nilpotentes dans~$\Lg^*$. Ce sont les \emph{paramètres d'Arthur unipotents} géométriques. Le mot \enquote{généralisé} signifie que l'orbite n'est pas forcément la régulière. Voici la construction de~$M$ ainsi que des données supplémentaires dans ces cas. Soit $O$ une telle orbite. Prenons un élément $\phi$ dedans. La forme de Killing induit $\Lg^*\ito\Lg$. Notons $f\in\Lg$ l'image de $\phi$. Selon Jacobson-Morozov, cet élément fait partie d'un $\frsl_2$-triplet $h,e,f$, unique à conjugaison près. Notons $\sigma:\SL_2\to\LG$ l'homomorphisme associé. Le diagonal $\GG_m$ induit une graduation sur $\Lg$. Nous savons que $\phi([\cdot,\cdot])$ est une forme symplectique sur sa partie de poids~1. Choisissons une lagrangienne dedans. On la combine avec la partie de poids au moins~2 de $\Lg$ pour former une algèbre de Lie nilpotente, notée~$\frus$. Notons par $\Us$ le groupe correspondant. La projection, notée $\bphi$, de $\phi$ dans $\frus^*$ est alors invariante sous $\Us$. La multiplication d'un côté de $\LG$ sur lui-même induit une action sur son espace cotangent. Il s'en suit l'application moment $\LG$-équivariante $\rmT^*\LG\to\Lg^*$. En particulier, elle est équivariante sous le diagonal $\Gm\subset\SL_2\overset{\sigma}{\to}\LG$, ainsi que sous le centralisateur $\Zs\subset\LG$ de $\sigma$. On va oublier cette action de $\LG$ en ne regardant que l'effet de son sous-groupe $\Us$.  On obtient ainsi la nouvelle application moment $\rmT^*\LG\to\frus^*$ qui hérite la $\Gm\times\Zs$-équivariance. Ensuite, on modifie l'action de $\Gm$ sur $\rmT^*\LG$ en la composant avec l'action scalaire au carré sur les vecteurs cotangents. Pour garder la $\Gm$-équivariance de l'application moment, on augmente de poids~2 son action sur $\frus^*$. Un bienfait est que $\bphi\in\frus^*$, initialement de poids~$-2$, devient point fixe sous cette nouvelle action de $\Gm$. Il faut rappeler maintenant qu'on dispose d'une deuxième application moment $\rmT^*\LG\to\Lg^*$ venant de l'action de l'autre côté de $\LG$. Cette application respecte l'action peu habituelle de $\Gm$ sur $\rmT^*\LG$ qu'on vient de décrire et son action conventionnelle sur $\Lg^*$ par scalaire au carré. Elle respecte aussi l'action ci-dessus de $\Zs$ sur $\rmT^*\LG$ et son action triviale sur $\Lg^*$. La variété symplectique, notée~$\Ms$, est prise comme la réduction hamiltonnienne de $\rmT^*\LG$  par rapport à $\Us$ au point $\bphi\in\frus^*$. Elle hérite une action de $\LG\times\Gm\:(\times\Zs)$ ainsi qu'une application moment vers $\Lg^*$ vérifiant toutes les conditions requises. Mentionnons en passant qu'elle joue un rôle central déjà dans la correspondance de Deligne-Langlands entre les représentations du groupe $\Zs$ et celles de l'algèbre classique de Hecke-Iwahori mixte. En effet, cette algèbre se réalise comme le K-groupe de $\Coh\left((\tcN\times_{\Lg^*}\tcN)/(\LG\times\Gm)\right)$ où $\tcN$ est la résolution de Springer du cône nilpotent $\cN\subset\Lg^*$. Par convolution, cette catégorie agit sur \[\Coh\left((\tcN\times_{\Lg^*}\Ms)/(\LG\times\Gm\times\Zs)\right)\] de façon à commuter avec celle de $\Rep(\Zs)$. Le K-groupe associé réalise alors ladite correspondance. Le L-faisceau qui s'en déduit est le faisceau de Whittaker spectral généralisé associé à~$O$, ou plutôt à~$\sigma$. On le note~$\Ls$, qui, sous l'action de $\Zs$, se décompose ensuite en la somme directe des $\chi^*\otimes\Lsk$, où $\chi$ parcourt les représentations irréductibles de~$\Zs$. 

Pour en parler de la structure de factorisation, il faut d'abord les transporter, via la dualité de Koszul, vers la catégorie de Hecke spectrale sphérique $\SphsLG$ des faisceaux ind-cohérents sur $\LS{D}\times_{\LS{\Dc}}\LS{D}$ où $D$ est le disque formel $\Spf(\CC[[t]])$, $\Dc=\Spec(\CC((t)))$ sa version avec l'origine ôté, et $\LS{D}$ désigne le champ de modules des $\LG$-systèmes locaux sur~$D$, etc. Celle-ci admet une version factorisable, simplement en remplaçant le disque formel --- et similairement pour le disque formel percé --- par le voisinage formel de l'union de plusieurs points se déplaçant librement sur une courbe lisse quelconque. La version Betti de cette transportation est expliquée dans~\cite{ben-zviRelativeLanglandsDuality2024}. Auparavant, Vincent Lafforgue l'a aussi fait dans un article non-publié. Nous l'adaptons au cadre de Rham pour convenir à notre but.

On doit réinterpréter la variété symplectique $\Ms$ définissant notre L-faisceau comme la réduction hamiltonnienne de $\rmT^*\left((\LG\times\Ga)/\Us\right)$ par rapport à $\Ga$ au point $1\in\frga^*$ où la deuxième composante de l'inclusion $\Us\subset\LG\times\Ga$ est induite par $\bphi\in\frus^*$. De cette façon, on voit que la catégorie $\Qcoh(\Ms/\LG)$ est équivalente à \[\cO_1\tMod\biggl(\Qcoh\Bigl((\LG\times\Ga)\backslash\rmT^*\left((\LG\times\Ga)/\Us\right)\Bigr)\biggr)\] où $\cO_1\in\Qcoh(\frga^*/\Ga)$ est le faisceau gratte-ciel au point 1 muni de sa structure $\Ga$-équivariante naturelle. Cette dernière catégorie s'écrit plus succinctement comme $\cO_1\tMod\bigl(\Qcoh([\Lg^*\times\frga^*\to\frus^*]/\Us)\bigr)$, où les termes du complexe entre crochets sont placés en degrés 0 et~1. Ces équivalences respectent la structure monoïdale. Remarquons que l'unité de la dernière catégorie est la tensorisation avec $\cO_1$ du faisceau structural du champ classique $[\Lg^*\times\frga^*\to\frus^*]/\Us$. La tonte sur ces catégories se décompose en deux étapes, détaillées dans les deux paragraphes suivants.

La première vient de l'action scalaire au carré sur la direction cotangente. En gros, son effet consiste, compte tenu de la dualité de Koszul, à remplacer les espaces cotangents genre $V^*$ par les espaces tangents décalés~$V[-1]$. En même temps, pour toute application $V_1\to V_2$ entres les espaces tangents, elle remplace l'action naturelle de 
$\Qcoh(V_1^*)$ sur $\Qcoh(V_2^*)$ par l'action de $\Indcoh(V_1[-1])$ sur $\Indcoh(V_2[-1])$ induite par la convolution $V_1[-1]\times V_2[-1]\to V_2[-1]$.
Par exemple, $\Qcoh(\frga^*/\Ga)$ se transforme en la catégorie de Hecke spectrale pour $\Ga$ $\HksGa:=\Indcoh(\frga[-1]/\Ga)$. 
De même, elle transforme $\Qcoh(\Lg^*/\LG)$ en $\SphsLG$, et transforme $\Qcoh([\Lg^*\times\frga^*\to\frus^*]/\Us)$ en la catégorie $\Indcoh$ sur le champ dérivé 
$[\frus\to\Lg\times\frga]/\Us$, les termes de ces deux complexes étant  placés en degrés 0 et~1. Ce dernier champ s'interprète directement comme le produit fibré \[\LSUs:=\LS[\LG\times\Ga]{D}\times_{\LS[\LG\times\Ga]{\Dc}}\LS[\Us]{\Dc}.\] La catégorie $\Indcoh$ là-dessus acquiert ainsi une structure factorisable et admet l'action de $\SphsLG\times\HksGa$. Bien que la dualité de Koszul respecte les structures monoïdales des deux côtés, leur nature sont différentes, à tel point que l'unité de la catégorie $\Indcoh$ sur ce dernier champ est $\pi_*\omega$ où \[\pi:\LS[\Us]{D}\to\LSUs\] est l'application évidente et $\omega$ est le faisceau dualisant sur la source. D'autre part, faites attention que le faisceau $\cO_1$ n'est pas encore $\Gm$-équivariant puisque le point $1\in\frga^*$ bouge sous cette action. Il faut attendre la deuxième action.

La deuxième étape provient de l'action partout par conjugaison du diagonal $\Gm\subset\SL_2\overset{\sigma}{\to}\LG$. Il agit notamment de poids 2 sur notre $\Ga$, qui, rappelons-le, est un quotient de $\Us$. Son action sur $\frga^*$ est alors de poids~$-2$. Cette action transforme la catégorie $\HksGa$ en sa version tondue notée $\Hkst$ comme d'habitude. Cette tonte n'affecte pas $\SphsLG$ puisque $\Gm$ agit sur $\LG$ par conjugaison intérieure, si bien que ses actions sur les divers champs des $\LG$-fibrés principaux sont toutes équivalentes à l'action triviale.

Ces deux actions commutent entre elles. Sous leur composée, le faisceau $\cO_1$ devient $\Gm$-équivariant, ce qui nous permet de parler de sa version tondue, notée $\cE$, dans $\Hkst$. On sait que l'objet $\cE$ est naturellement factorisable.

On est prêt à donner la forme factorisable sur $\SphsLG$ de notre L-faisceau~$\Ls$. Rappelons qu'initialement, on l'a définit dans $\SatLG$ comme image directe du faisceau structural tondu $\OMst$, l'unité de $\Qcoh(\Ms/\LG)\tonte$. En regardant cette dernière catégorie comme module sous l'action par tensorisation de la catégorie monoïdale $\SatLG$, on peut interpréter le L-faisceau comme l'endomorphisme interne reletif de $\OMst$. D'après les analyses qu'on vient d'effectuer sur l'effet de la dualité de Koszul, $\Qcoh(\Ms/\LG)\tonte$ est équivalente à \[\cE\tMod\Bigl(\Indcoh\bigl(\LSUs\bigr)\tonte\Bigr),\] de sorte que l'unité $\OMst$ correspond à $\cE\star(\pi_*\omegat)$, l'unité de l'autre côté, où~$\star$ désigne la convolution par $\Hkst$; en même temps, l'action par convolution de $\SphsLG$ sur cette nouvelle catégorie remplace celle par tensorisation de $\SatLG$ sur $\Qcoh(\Ms/\LG)\tonte$, si bien que le L-faisceau correspond à l'endomorphisme interne relatif à $\SphsLG$ de $\cE\star(\pi_*\omegat)$.  De cette façon, l'avatar du L-faisceau dans $\SphsLG$ acquiert une structure factorisable. Il y a plus, même si l'on n'en a pas besoin dans cet article: comme il est un objet associatif dans la catégorie des objets factorisables, il admet en fait une structure~$\rmE_3$. À l'aide du pont $\Ga\leftarrow\Us\to\LG$, on peut l'exprimer aussi comme l'image de~$\cE$ sous la composée $\Hkst\to\Hkst[\Us]\to\SphsLG$ de l'image réciproque et l'image directe. Cette réécriture permet de voir la fonctorialité du L-faisceau. Pour $\sigma:\SL_2\to\LG_1\to\LG_2$ où la dernière flèche modulo les centres est un isomorphisme, le L-faisceau associé à $\sigma$ pour~$G_2$ est l'image directe de celui pour~$G_1$ via $\SphsLG[\LG_1]\to\SphsLG[\LG_2]$, simplement parce que les $\Us$ correspondants sont isomorphes. Tout ceci reste valable pour la version factorisable.

Comme promis dans l'introduction, rendons explicite la structure de fusion entre ses composantes isotypiques~$\Lsk$. Notons la représentation triviale de~$\Zs$ par~1.

\begin{thm*}
L'objet $\Lsk[1]$ admet une structure factorisable sur laquelle les $\Lsk$ admettent tous une structure de module factorisable. Il existe un homomorphisme naturel de $\Hom_{\Zs}(\chi_1\otimes\chi_2,\chi_3)$ vers l'espace des opérateurs d'entrelacement de type $\left(\begin{array}{c}
          \Lsk[\chi_3]\\\Lsk[\chi_1],\Lsk[\chi_2]
    \end{array}\right)$.
\end{thm*}

En-voici la construction. Soit $C$ une courbe lisse quelconque. La formulation de toutes les structures factorisables est basée sur~$\Ran$, le préchamp dont la valeur sur un schéma affine de test $S$ est l'ensemble des sous-ensembles finis non-vides de $\Hom(S,C)$. On a les morphismes naturels $C^J\to\Ran$ pour tout ensemble fini non-vide $J$, de sorte que $\Ran$ en est la colimite par rapport aux immersions diagonales $C^{J'}\hookrightarrow C^J$ induites par les applications surjectives $J\twoheadrightarrow J'$. Notons $\Rani$ le sous-préchamp de $\Ran^2$ classifiant deux sous-ensembles finis non-vides de $\Hom(S,C)$, le premier contenu dans le deuxième. Notons par $\pr_1: \Rani\to\Ran$ la projection au premier sous-ensemble et par $\pr_2$ celle au deuxième. Remarquons que $\pr_2$ se décompose en l'inclusion $\Rani\hookrightarrow\Ran^2$ suivie de l'union vers $\Ran$.

On a besoin de la notion de catégorie factorisable commutative comme décrite dans la note intitulée \enquote{Chiral Categories} de Sam Raskin. Une telle catégorie $\cCRan$ est munie des applications supplémentaires $\cCRan\boxtimes\cdots\boxtimes\cCRan\to\cup^*\cCRan$ par rapport à l'union $\cup:\Ran\times\cdots\times\Ran\to\Ran$, étendant l'équivalence structurale sur $(\Ran\times\cdots\times\Ran)_\disj$. Pour toute catégorie monoïdale symétrique~$\cC$, on peut construire canoniquement une telle catégorie dont la restriction à $C\dR$ est $\Dmod(C)\otimes\cC$.
On sait que $\Indcoh(\LS[\Zs]{\DRan})$ est une telle catégorie factorisable commutative associée à la catégorie symétrique $\Rep(\Zs)$.
Pour tout $J$ fini non-vide, prenons $\Ran_J:=C^J\dR\times_{\Ran,\pr_1}\Rani$, considéré comme au-dessus de $\Ran$ via~$\pr_2$. C'est le préchamp de base entrant dans la formulation de la notion de module factorisable. Notons $\Ran_{J,\disj}$ le sous-préchamp de $\Ran_J$ où les diviseurs indexés par $J$ sont disjoints. On remplacera systématiquement $\Ran$ par $\Ran_J$ dans les notations pour désigner leur changement de base par rapport à $\Ran_J\to\Ran$; de même lorsqu'on le remplace par $\Ran_{J,\disj}$, etc. La propriété commutative qu'on vient d'évoquer entraîne alors une action de $\Rep(\Zs)^{\otimes J}$ sur $\Indcoh(\LS[\Zs]{\DRan[J]})$. 

Il est clair que pour analyser les composantes du L-faisceau, il faut prendre l'endomorphisme interne de l'unité de $\Qcoh\left(\Ms/(\LG\times\Zs)\right)\tonte$ relatif à l'action de $\SatLG\otimes\Rep(\Zs)$.
Comme précédemment, sa version factorisable $\LsZRan$ dans $\SphsLGRan\times_{\Ran}\LS[\Zs]{\DRan}$ s'obtient grâce à la dualité de Koszul, comme l'endomorphisme interne relatif de l'unité $\cERan\star(\pi'_*\omega^{',\fatslash})$ de \[\cERan\tMod\Bigl(\Indcoh\bigl(\LSZUsRan\bigr)\tonte\Bigr),\] où $\pi':\LS[\Zs\Us]{\DRan}\to\LSZUsRan$ est le morphisme naturel et $\omega'$ est le faisceau dualisant sur la source.  À l'aide cette fois du pont $\Ga\leftarrow\Zs\Us\to\Zs\times\LG$, on peut l'exprimer aussi comme l'image de $\cERan$ sous \[\HkstRan\to\SphstZRan{\Zs\Us}\to\SphsZRan{\Zs\times\LG}\to\SphsLGRan\times_{\Ran}\LS[\Zs]{\DRan}.\] Son image directe, notée $\LskRan{1}$, dans $\SphsLG$ est la version factorisable de $\Lsk[1]$. D'autre part, sous l'action de $\Rep(\Zs)^{\otimes J}$, l'objet  $\LsZRan[J]$ dans $\SphsLGRan[J]\times_{\Ran_J}\LS[\Zs]{\DRan[J]}$ engendre, pour toute famille de représentations irréductibles $(\chi_j)_{j\in J}$ de~$\Zs$, l'objet $(\boxtimes\chi_j)*\LsZRan[J]$, dont l'image directe dans $\SphsLGRan[J]$ sera notée par $\LskRan[J]{\boxtimes\chi_j}$, laquelle fournit la structure de fusion qu'on veut. Pour $J=\{*\}$ un singleton, l'objet $\LskRan[*]{\chi}$ munit $\Lsk[\chi]$ d'une structure de module factorisable sur $\Lsk[1]$. Pour $J=\{1,2\}$ une paire, considérons l'objet $\LskRan[1,2]{\chi_1\boxtimes\chi_2}$. On voit que sa restriction à l'ouvert où les deux diviseurs indexés sont disjoints  est la fusion externe de $\LskRan[*]{\chi_1}$ et $\LskRan[*]{\chi_2}$ d'après un argument hyper-descente, et que sa restriction au diagonal où les deux éléments indexés de $\Hom(S,C)$ coïncident est \enquote{$\LskRan[*]{\chi_1\otimes\chi_2}$} qui se décompose en la somme directe des $\Hom_{\Zs}(\chi_3,\chi_1\otimes\chi_2)\bigotimes\LskRan[*]{\chi_3}$ où $\chi_3$ parcourt les représentations irréductibles de~$\Zs$. On met les guillemets pour signaler que $\chi_1\otimes\chi_2$ n'est pas irréductible en générale. On a aussi la fonctorialité. Pour $\sigma:\SL_2\to\LG_1\to\LG_2$ où la dernière flèche modulo les centres est un isomorphisme, on note les centralisateurs $\Zs^1,\Zs^2$. Alors pour toute famille de représentations irréductibles $\chi^1_j$ de $\Zs^1$, l'image directe sous $\SphsZRanJ[J]{\LG_1}\to\SphsZRanJ[J]{\LG_2}$ de $\LskRan[J]{\boxtimes\chi^1_j}$ est la somme directe des $\LskRan[J]{\boxtimes\chi^2_j}$ où $\chi^2_j$ parcourt les représentations irréductibles de $\Zs^2$ dont la restriction à $\Zs^1$ est~$\chi^1_j$.

Il est temps de passer à la grassmannienne affine $\Gr_G:=G(\Dc)/G(D)$. Grâce à
l'équivalence de Satake dérivée, voir~\cite{bezrukavnikovEquivariantSatakeCategory2008}, la catégorie de Hecke spectrale sphérique s'identifie avec la catégorie de Hecke sphérique $\SphG$ des D-modules sur le quotient $G(D)\backslash G(\Dc)/G(D)$. On note $\Fs,\Fsk$ les objets dans $\SphG$  correspondant à $\Ls,\Lsk$. On a en fait aussi les objets factorisables $\FsRan, \FskRan{1}$ dans la catégorie $\SphGRan$ des D-modules sur $G(\DRan)\backslash G(\DcRan)/G(\DRan)$, parce que cette équivalence admet elle même une version factorisable récemment montrée dans~\cite{campbellLanglandsDualityBeilinsonDrinfeld2024}. Le D-module $\Fs$ est fonctoriel au sens suivant. Pour $\sigma:\SL_2\to G_2\to G_1$ où la dernière flèche modulo les centres est un isomorphisme, l'image réciproque de \enquote{$\Fs$ pour $G_1$} est bien \enquote{$\Fs$ pour $G_2$}. La raison est que dans cette situation, l'image réciproque correspond à l'image directe au côté spectral, et on conclut par la fonctorialité de~$\Ls$ qu'on a montrée. De même, on a $\FskRan[J]{\boxtimes\chi_j}$ dans $\SphGRan[J]$, de sorte que $\FskRan[*]{\chi}$ munit $\Fsk$ d'une structure de module factorisable sur $\Fsk[1]$. Suivant la fonctorialité des $\LskRan[J]{\boxtimes\chi_j}$ qu'on vient de décrire, les D-modules $\FskRan[J]{\boxtimes\chi_j}$ en jouissent aussi sous les images réciproques.

On peut maintenant donner la construction des algèbres vertex associées à~$\sigma$. On emploie le pluriel ici parce qu'il faut encore la donnée d'un fibré en droites factorisables sur $\GrRan$. L'algèbre vertex $\VsRan$ est alors simplement la section globale relative sur $\GrRan$, c'est-à-dire l'image directe sous $\GrRan\to\Ran$, de $\FsRan$ tordu par ce fibré. Énumérons rapidement ses propriétés. Sa version sur un point, notée $\Vs$, se décompose sous l'action de $\Zs$ comme la somme directe des $\chi^*\otimes \Vsk$. Ici, $\Vsk[1]$ est une sous-algèbre vertex puisqu'elle admet une version factorisable $\VskRan{1}$. De plus, les composantes $\Vsk$ admettent toutes une structure de module sur~$\Vsk[1]$ et fusionnent de la même manière que les représentations de~$\Zs$, grâce aux objets $\VskRan[J]{\boxtimes\chi_j}$ sur $\Ran_J$.

Un petit souci: il n'y a pas assez de fibrés en droites factorisables sur $\GrRan$. Autrement dit, on va rater beaucoup de niveaux. Pour les récupérer, au lieu d'utiliser ces fibrés sur $\GrRan$, on essaie d'utiliser tous ceux sur $\GrRan[\Gsc]$ où $\Gsc$ est le recouvrement simplement connexe du groupe semi-simple $[G,G]$ des commutateurs. Ceci est possible grâce au fait que $\Fsk$ est supporté sur une seule composante connexe de $\Gr_G$ comme on l'a mentionné déjà dans l'introduction. Plus précisément, les composantes de $\Gr_G$ sont naturellement indexées par les caractères du centre de~$\LG$. Ce centre étant contenu dans~$\Zs$, toute~$\chi$ irréductible induit un tel caractère par restriction. On peut alors parler de la $\chi$-composante connexe. Le D-module $\Fsk$ est supporté sur cette composante, ce qui permet de situer le support de $\FskRan[J]{\boxtimes\chi_j}$. Ensuite, prenons un tore maximal~$T$ de~$G$. L'action par conjugaison de ce tore sur $G$ induit canoniquement son action sur~$\Gsc$. Il s'en suit un homomorphisme $\Gsc\rtimes T\to G$, lequel implique une application surjective entre les grassmanniennes affines des deux côtés.  Pour tout $j\in J$, on a une section canonique pour $\GrRanJ[J]{\Gm}\to\Ran_J$, donnée de la façon suivante. Soit $S$ un schéma de test et $S\to\Ran_J$ une application. En particulier on a des applications $S\to C$ indexées par~$J$, lesquelles induisent des diviseurs sur $C\times S$ notés~$\Gamma_j, j\in J$.  Le fibré en droites $O(\Gamma_j)$ sur $C\times S$, banalement trivialisé en dehors de tous les diviseurs associés à l'application $S\to\Ran_J$, fournit cette section canonique. Pour tout homomorphisme $\lambda_j:\Gm\to T$, on en déduit une section de $\GrRanJ[J]{T}\to\Ran_J$ par extension de groupes. On note $\tlambdaj:\Ran_J\to\GrRanJ[J]{T}$ le produit de ces sections indexées par~$J$. On obtient le diagramme 
\[\begin{array}{ccccc}
     \Grlambda &\to &\GrRanJ[J]{\Gsc\rtimes T} &\to &\GrRanJ[J]{G}  \\\downarrow&&\downarrow &&\\
     \Ran_J&\to &\GrRanJ[J]{T} && 
\end{array}
\] où le carré est le changement de base par rapport~$\tlambdaj$. On sait que $\Grlambda$ est canoniquement isomorphe à $\Gsc(\DcRan[J])/\Ad(\tlambdaj)\Gsc(\DRan[J])$. Considérons l'application composée de la première ligne. Remarquons que son image contient le support de $\FskRan[J]{\boxtimes\chi_j}$ lorsque pour tout~$j$, $\lambda_j$ indexe aussi la $\chi_j$-composante de~$\Gr_G$. 

D'autre part, on sait que  tout fibré en droites factorisable sur $\GrRan[\Gsc]$ s'étend de façon unique à un fibré en droites factorisable sur $\GrRan[\Gsc\rtimes T]$ dont la restriction sur $\GrRan[T]$ est trivialisée. Par restriction, cela donne un fibré en droites sur $\Grlambda$. 

On va redéfinir $\VskRan[J]{\boxtimes\chi_j}$ comme l'image directe sous $\Grlambda\to\Ran_J$ de l'image réciproque de $\FskRan[J]{\boxtimes\chi_j}$ tordue par ce fibré. Ceci est indépendant du choix des $\lambda_j$ compatibles aux $\chi_j$ au sens ci-dessus. Ces objets nous redonne l'algèbre vertex $\Vsk[1]$ ainsi que ses modules $\Vsk$ vérifiant les mêmes propriétés, mais cette fois associés à tout fibré en droites factorisable sur $\GrRan[\Gsc]$ comme on désire. 

Faites attention que lorsque le fibré sur $\GrRan[\Gsc]$ vient d'un, noté~$L$, sur $\GrRan$ par restriction, cette nouvelle définition n'est compatible avec l'ancienne que pour~$G$ semi-simple. Sinon le réseau des copoids centraux de~$G$ donnera une algèbre vertex de type Heisenberg qui sera perdue en passant de $\Gr_G$ à~$\Gr_{\Gsc}$. Pour la récupérer, il faudrait remplacer $\Gsc$ par le recouvrement universel de~$G$, ou plutôt par ses approximations $Z'\times\Gsc$ où $Z'$ est un recouvrement de la composante neutre du centre de~$G$ puisque le recouvrement universel de~$G$ en est la limite. Pour~$G$ semi-simple, les deux versions de $\VskRan[J]{\boxtimes\chi_j}$ sur $\Ran_J$ diffèrent d'un fibré en droites qui est la restriction de~$L$ via $\Ran_J\to\GrRanJ[J]{T}\to\GrRanJ[J]{G}\to\GrRan$ où la première flèche est la section~$\tlambdaj$. Cette différence n'a pas d'importance. 

En tout cas, on supposera dans le suivant que~$G$ est semi-simple. Notre but dans les paragraphes prochains est de donner des formules explicites aux $\Vsk$ dans des cas favorables. Pour l'achever, il suffit de suppose que notre courbe~$C$ est~$\PP^1$, et que~$G$ est presque simple et même de type adjoint pour que $\Gr_G$ aie le plus de composantes. On les suppose. Alors, d'après la conjecture 3.3.7 dans~\cite{gaitsgoryParametersDualityMetaplectic2018} montrée dans~\cite{taoCentralExtensions2021}, les fibrés en droites factorisable sur $\GrRan[\Gsc]$ sont classifiés par le niveau, un entier, noté~$\kappa$. On fixe un Borel $B\subset G$ contenant le tore maximal~$T$.

Comme expliqué au début, au lieu de travailler directement sur les D-modules sous-jacents sur la grassmannienne affine, on prend d'abord les D-modules automorphes qui sont leurs transformés de Radon longue. C'est légitime car ils donnent les mêmes sections globales tordues. Ces D-modules automorphes correspondent aux faisceaux spectraux sous Langlands géométrique. Ce sont des D-modules sur $\Fib_G$, le champ des modules des $G$-fibrés principaux sur la courbe~$C$. D'après le  principe tracé dans~\cite{ben-zviRelativeLanglandsDuality2024}, ces D-modules auraient dû être des \enquote{faisceaux de périodes}. Or, la recette existante pour les décrire demande que la variété symplectique définissant le L-faisceau soit hyper-sphérique, ce qui n'est pas notre cas sauf si l'orbite~$O$ est la régulière parce que la $\LG$-variété $\LG/\Us$ n'est pas sphérique pour les autres orbites.

Lorsque l'orbite est la régulière, on sait que le D-module automorphe associé est le faisceau constant sur $\Fib_G$. Notons $\Fibx$ le champ des $G$-fibrés principaux sur la courbe munis d'une trivialisation sur le voisinage formel $D_x$ de~$x$. Il est isomorphe au quotient $G(\Dc_x)/G(C-x)$. Le centralisateur $\Zs$ est égal au centre de~$\LG$.  Ses caractères sont en bijection à la fois avec les composantes de $\Gr_G$ et avec celles de $\Fibx$. On note $\lambda$ le copoids dominant minuscule associé au caractère~$\chi$. On sait qu'il détermine un unique élément $\epsilon$ du groupe de Weyl de $\Gsc$ tel que l'action par conjugaison de $\gamma:=t^{-\lambda} \epsilon$ préserve le \enquote{Borel standard} de l'extension centrale minimale du groupe de lacets $\Gsc(\Dc)$ correspondant au fibré en droites de niveau~1 sur $\Gr_{\Gsc}$. L'élément $\epsilon$ préserve la réunion des racines simples et de $-\theta$, où~$\theta$ est la plus haute racine. Il est caractérisé comme celui qui envoie $-\theta$ à la racine simple correspondant à~$\lambda$. L'action de~$\gamma$ permute les coracines simples affines ainsi que les poids fondamentaux affines. On fixe une identification de $D_x$ avec~$D$ respectant les origines.
La $\chi$-composante de $\Gr_G$, resp. $\Fibx$,  est isomorphe à $\Gsc(\Dc)/\Ad(\gamma)\bigl(\Gsc(D)\bigr)$, resp. $\Gsc(\Dc_x)/\Ad(\gamma)\bigl(\Gsc(C-x)\bigr)$. Notre fibré en droites descend sur ces composantes. Supposons $\kappa$ positif, alors $\Vsk$ est la section globale dérivée de ce fibré sur la $\chi$-composante de $\Fibx$. On en déduit que $\Vsk[1]$ est~$\LkL{0}$, le quotient simple du vacuum de niveau~$\kappa$. On le note plus simplement par~$\Lk$. 
L'action de $\gamma$ envoie le plus haut poids $\kappa\Lhat_0$ du vacuum à $\kappa\Lhat_i$ où $\Lhat_i$ est le poids fondamental affine dont la restriction, qui est $\Lambda_i$ pour $\gamma$ non-trivial et zéro sinon, au tore maximal de $\Gsc$ est le poids fondamental cominuscule correspondant à~$\lambda$. Alors $\Vsk$ est $L(\kappa\Lhat_i)$, le module simple de l'algèbre de Lie affine de plus haut poids $\kappa\Lhat_i$. Pour~$\chi$ non-trivial, il est $\LkL{\Lambda_i}$, le quotient simple du module de Weyl affine associé à la représentation irréductible cominuscule de $\Gsc$ de plus haut poids~$\Lambda_i$. Comme la fusion de $\Vsk[\chi_1]$ et $\Vsk[\chi_2]$ est $\Vsk[\chi_1\chi_2]$, ces modules de l'algèbre vertex $\Lk$ sont tous des courants simples, c'est-à-dire des modules inversibles.

En général, nous allons montrer qu'un lien similaire, quoique moins précis, demeure entre ces D-modules automorphes et les faisceaux pervers simples. En fait, pour ce qui nous concerne, on considérerait plutôt leurs images réciproques sur $\FibxU$, le champ des $G$-fibrés principaux sur~$C$ munis d'une réduction à~$U$ à un point $x\in C$ fixé où~$U$ est le radical unipotent du Borel~$B$. On choisit de travailler sur ce nouveau terrain parce qu'on aura recours aux techniques qui ne sont disponibles que pour une catégorie de plus haut poids, et on sait que $\Dmod(\FibxU)$ en est un exemple: en effet, c'est une version affine de la catégorie des D-modules Borel-monodromiques sur une variété de drapeaux partiels. On note $\As$ le D-module monodromique correspondant à~$\Ls$.

Pour détecter les sous-quotients simples de ces D-modules monodromiques, on doit calculer $\Hom(\Psph,\cdot)$ où $\Psph$ est la couverture projective d'un objet simple. L'exposant ici indique qu'elle est déjà \enquote{sphérique à droite}.
On a besoin du diagramme commutatif suivant où toutes les flèches verticales sont des équivalences:
\[
\begin{tikzcd}[column sep=small]
\FibxU\ar[rr]\ar[d,"\text{Radon}"]&&\Dmod(\Fib_G)\ar[d,"\text{Radon}"]\\
\Dmod(I_u\backslash \Gr_G)\ar[d,"\text{Koszul}"]\ar[rr]&&\SphG\ar[dd,"\text{Satake}"]\\\Dmod\bigl(I\backslash G(\Dc)/(I_u^-,\Psi)\bigr)\tonte\ar[d,"\text{\cite{arkhipovPerverseSheavesAffine2009}}"]& \\\Qcoh(\tcN/\LG)\tonte\ar[r]&\Indcoh(\cN/\LG)\tonte\ar[r]&\Qcoh(\Lg^*/\LG)\tonte.
\end{tikzcd}
\]
Ici $I\subset G(D)$ est le sous-groupe d'Iwahori prenant valeur dans~$B$ à l'origine avec radical pro-unipotent~$I_u$. Similairement, $I^-\subset G(D)$ est le sous-groupe prenant valeur dans le Borel opposé $B^-$ à l'origine avec radical pro-unipotent~$I_u^-$. Comme d'habitude, $\tcN\to\cN$ désigne la résolution de Springer. Le caractère additif $\Psi:I_u^-\to\Ga$ est induit d'un caractère additif générique du radical unipotent de~$B^-$. La notation $(I_u^-,\Psi)$ sert à exprimer la condition d'Iwahori-Whittaker, c'est-à-dire la $\Psi$-équivariance sous l'action de~$I_u^-$. La dualité de Koszul ici est celle proposée dans~\cite{beilinsonWallcrossingFunctorsModules1999} et montrée dans~\cite{bezrukavnikovKoszulDualityKacMoody2013}. Elle est en réalité une équivalence entre les versions mixtes des deux côtés, sous laquelle les torsions à la Tate des deux côtés ne coïncident pas mais diffèrent par un décalage cohomologique. En oubliant ces torsions, on arrive à l'équivalence qu'on emploie ici où la tonte apparaît.  

Revenons à notre objet projectif~$\Psph$. D'après~\cite{beilinsonTiltingExercises2004}, il se transforme sous Radon en un faisceau basculant indécomposable qui, après l'effet de la dualité de Koszul, devient la tonte d'un faisceau pervers d'Iwahori-Whittaker simple. Ce faisceau simple, sous l'équivalence d'Arkhipov-Bezrukavnikov qui est t-exacte par rapport à une t-structure \enquote{perverse} naturelle au côté cohérent, correspond à un objet simple dans le cœur sur $\tcN/\LG$. Sa projection sur $\cN/\LG$ est nulle si le quotient simple, qui est déjà sphérique à droite, de $\Psph$ n'est pas en même temps \enquote{sphérique à gauche}, autrement dit si ce quotient ne vient pas de $\Fib_G$. Sinon, cette projection, notée~$S$, reste un objet simple pour une t-structure perverse cohérente naturelle et on a, d'après le diagramme, que
\[\Hom(\Psph,\As)\ito\Hom(S\tonte,\Ls).\]

Le support d'un tel objet simple sur $\cN$ est irréductible; l'objet lui-même est déterminé par sa restriction sur l'orbite maximale dans son support. Cette restriction se révèle non-dérivée et irréductible au sens habituel. De cette façon, les D-modules monodromiques simples et sphériques sont en fait en bijection avec les couples $(O',\chi')$ où~$O'$ est une orbite nilpotente dans $\Lg^*$ et~$\chi'$ un faisceau cohérent équivariant irréductible sur~$O'$. On note le D-module $\IC_{O',\chi'}$. C'est une forme de la bijection de Lusztig-Vogan.

L'application moment définissant notre L-faisceau s'interprète plus explicitement comme l'application équivariante évidente $\LG\times\Sls\to\Lg^*$. Ici $\Sls:=phi+\ker(\ad e)\subset\Lg^*$ est une tranche de Slodowy transversale à~$O$. Ceci implique que toutes les orbites nilpotentes apparaissant dans son image contiennent $O$ dans leur adhérence. On doit donc avoir $O\subseteq\bar{O'}$ pour que le $\Hom$ ci-dessus soit non-nul. On dit~$O'$ est plus grande que~$O$ et écrit $O\leq O'$.
 
Lorsque $O=O'$, le $\Hom$ se calcule sans beaucoup de peines puisque le changement de base de l'application moment par rapport à l'inclusion $O\hookrightarrow\Lg^*$ est juste la fibration $\LG\to O$ envoyant l'origine à~$\phi$. On sait que le centralisateur de $\phi$ est le produit semi-direct de son radical unipotent et le groupe réductif éventuellement non-connexe $\Zs$. Un faisceau cohérent équivariant irréductible $\chi$ sur $O$ correspond alors naturellement à une représentation irréductible de~$\Zs$. On notera cette représentation aussi par~$\chi$, par abus de notation. À un décalage ne dépendant que de l'orbite $O$ près, le $\Hom$ est alors isomorphe à~$\chi^*$, le dual de~$\chi$. D'autre part, on a vu que l'action de $\Zs$ induit la décomposition de $\Ls$ en la somme directe des $\chi^*\otimes\Lsk$ où $\chi$ parcourt les représentations irréductibles de~$\Zs$, de sorte que $\Hom(S\tonte,\Lsk)$ est isomorphe à~$\CC$ (décalé) si $S$ correspond à $(O,\chi)$, et nul si $S$ correspond à d'autres faisceaux équivariants irréductibles sur~$O$. Ensuite, en transposant au côté automorphe l'action de $\Zs$, on obtient une décomposition de $\As$ en la somme directe des $\chi^*\otimes\Ask$ où $\Ask$ est isomorphe à $\IC_{O,\chi}$ (décalé) modulo les objets simples correspondant aux orbites strictement plus grandes que~$O$.

Supposons cette fois $\kappa$ seulement strictement plus grand que critique. La description ci-dessus des $\Ask$ en termes des objets simples nous sert de levier pour donner des formules explicites aux $\Vsk$ qui en sont les sections globales $\kappa$-tordues sur $\Fibx$, lequel est  une version affine de la variété de drapeaux partiels. Fixons un deuxième point $y$ de notre courbe différent du point~$x$. On peut \enquote{compléter les drapeaux} en choisissant une $B^-$-réduction au point $y$ du $G$-fibré correspondant et obtenir un schéma au-dessus de $\Fibx$. On le note $\Fibhxy$. Il est isomorphe au quotient $G(\Dc_x)/\II^-$ où $\II^-$ est le sous groupe de $G(C-x)$ prenant valeur dans $B^-$ au point~$y$. Les sections globales tordues là-dessus de nos D-modules donnent les mêmes résultats que celles sur $\Fibx$ puisque la cohomologie cohérente dérivée de $G/B^-$ est~$\CC$. Cette réduction supplémentaire au point $y$ invite l'action de Hecke-Iwahori. On a aussi besoin du champ noté $\FibxBy$ qui classifie les $G$-fibrés sur la courbe munis à la fois d'une $B$-réduction au point $x$ et d'une $B^-$-réduction au point~$y$. Introduisons aussi le champ $\FibxUy$ dont le sens est prononcé déjà par la notation. Le diagramme commutatif
\[
\begin{tikzcd}
    \Fibhxy\arrow[r]\arrow[d] & \FibxUy\arrow[r]\ar[d] & \FibxBy\arrow[d]  \\  \Fibx\ar[r]&\FibxU\arrow[r]& \Fib_G
\end{tikzcd}
\]
résume la relation entre ces champs.

Rappelons la catégorie monoïdale de Hecke-Iwahori $\Hk_G$. Elle est constituée des D-modules sur $I\backslash G(\Dc)/I$. On fixe un automorphisme de $G$ échangeant $B,B^-$ et induisant l'automorphisme sur $T$ donné simplement par l'inverse. Il permet d'identifier $\Hk_G$ avec $\Dmod\left(I^-\backslash G(\Dc)/I^-\right)$ qui agit naturellement sur $\Dmod(\FibxUy)$ par convolution au point~$y$. Ici on a fixé une identification du voisinage formel $D_y$ là-dessus avec~$D$ respectant les origines. On prend la version à droite de l'action. Cette action fait partie du diagramme suivant:
\begingroup
\renewcommand{\arraystretch}{1.5}
\[
\begin{array}{ccl}
    \Dmod(\FibxUy) &\curvearrowleft&\:\:\Hk_G \\
    \Big\downarrow \text{\scriptsize Radon}& &\\
     \Dmod\bigl(I_u\backslash G(\Dc)/I\bigr)&\curvearrowleft&\:\:\Hk_G\\
\Big\downarrow\text{\scriptsize Koszul}&&\quad\Big\downarrow\text{\scriptsize Koszul}\\
  \Dmod\bigl(I\backslash G(\Dc)/I_u\bigr)\tonte   &\curvearrowleft&\HkGmonodrom\tonte.
\end{array}
\]
\endgroup
Ici toutes les flèches verticales sont des équivalences. Les lignes pointillées signifient \enquote{monodromique par rapport à l'action de $I$}, une condition plus faible que l'équivariance.

Le groupe de lacets $G(\Dc)$ est stratifié par $IwI$ où $w$ parcourt le groupe de Weyl affine étendu. Le faisceau constant là-dessus s'étend au D-module standard~$\Delta_w$, au D-module costandard~$\nabla_w$, ainsi qu'au D-module simple~$C_w$. Un fait important pour nous est que la transformée de Radon dans ce diagramme commute avec l'action de ces D-modules standards. On montre qu'elle commute avec les costandards qui sont les inverses des standards sous la convolution. Rappelons que cette transformation, ou plutôt son inverse qui va dans l'autre direction, vient de la convolution au point $x$ sur le D-module standard supporté sur la grande strate de $\FibxBy$, celle du $G$-fibré trivialisable sur la courbe $C$ muni des $B,B^-$-réductions, dont la position relative est générique, aux points~$x,y$. Il est bien connu que sous la convolution par $\nabla_w$ au point~$x$, ce D-module devient le D-module standard supporté sur une autre strate, à savoir la strate $w$ sous l'identification de $\FibxBy$ avec $I\backslash G(\Dc_x)/\II^-$. La convolution à droite par $\nabla_w$ au point~$y$ aura le même effet, ce qui montre ce qu'on veut.
 
Expliquons le rôle de cette action par les standards.
Notons $\rhohat$ la somme des poids fondamentaux affines. On dit qu'un poids affine (entier) est $\rhohat$-dominant si ses valeurs sur toutes les coracines simples affines sont au moins~$-1$, autrement dit, si sa somme avec $\rhohat$ est dominante. On a une meilleure connaissance du foncteur section globale tordue par un tel poids. Par exemple, le foncteur est t-exacte dans ce cas. Lorsque la valeur d'un poids $\rhohat$-dominant sur une coracine simple affine atteint le minimum~$-1$, on dit que c'est une coracine simple singulière pour ce poids. L'ensemble des réflexions simples singulières associées est noté $\Sing(\Lhat)$. La section globale tordue par un poids $\rhohat$-dominant $\Lhat$ sur $\Fibhxy$ annule le D-module simple $\IC_w$ sur l'adhérence de la strate $w$ de $\FibxBy$ si $w$ n'est pas de longueur maximale par rapport à la multiplication à droite par $\Sing(\Lhat)$; sinon elle donne le module simple $L(w\cdot\Lhat)$ où $w\cdot \Lhat:=w(\Lhat+\rhohat)-\rhohat$ est l'action décalée.

Le poids initial associé à notre niveau~$\kappa$ est $\kappa\Lhat_0$. Il n'est $\rhohat$-dominant que lorsque $\kappa\geqslant -1$. D'ailleurs, un tel poids $\rhohat$-dominant est singulier si et seulement si la cohomologie du fibré en droites associé sur $\Fibx$ est nulle. Ceci correspond à $\kappa=-1$. Dans ce cas, on peut déjà calculer $\Vsk$ pour $O$ sous-régulière car l'ambiguïté dans notre description du D-module automorphe sous-jacent $\Ask$ est constituée exactement des copies du faisceau constant sur la $\chi$-composante de $\Fibx$. On obtient  $\Vsk\ito L\bigl(w\cdot (-\Lhat_0)\bigr)$ où $w$
correspond au couple $(O,\chi)$. Cette formule ressemble tout à fait à celle pour l'orbite régulière qu'on a présentée.

En général, le niveau est seulement supposé strictement plus grand que critique. Il existe néanmoins un élément $v$ dans le groupe de Weyl affine tel que $v\cdot\kappa\Lhat_0$ est $\rhohat$-dominant. Ce nouveau poids est indépendant du choix de~$v$. On le notera par~$\kp$. On va travailler avec l'unique $v$ de longueur minimale.  D'après~\cite[\Thmdot{Cor}{1.6.2}]{kashiwaraKazhdanLusztigConjectureAffine1996}, le module $\Vsk$ est aussi la section globale $\kp$-tordue de $\Ask*\Deltavi$. La convolution ici est par rapport à l'Iwahori. Il faut signaler que comme Kashiwara et Tanisaki considèrent des D-modules à droite dans leur article, un niveau plus petit que critique pour nous paraîtra plus grand que critique là-bas. Comme $\Dmod(\FibxUy)$ est une catégorie de plus haut poids et un poids $\rhohat$-dominant interagit bien avec les objets simples dedans, on aimerait bien déterminer les sous-quotients simples de cette convolution. Pour cela, il faut calculer $\Hom(P_z,\Ask*\Deltavi)$ où $P_z$ est le recouvrement projectif d'un objet simple $\IC_z$. On ne s'intéresse qu'aux~$z$ de longueur maximale dans l'ensemble $z\Sing(\kp)$ puisque les autres ne contribuent pas dans la section globale tordue. On dit qu'un tel~$z$ \enquote{peut terminer par toutes les singularités de $\kp$}. Ce Hom est bien sûr pareil que $\Hom(P_z*\nabla_v,\Ask)$. Le projectif $P_z$ est transformé sous Radon en un objet basculant indécomposable qui se transforme ensuite en la tonte de l'objet simple $C_z$ sous la dualité de Koszul. Parallèlement, l'objet costandard $\nabla_v$ dans $\Hk_G$ devient la tonte de l'objet costandard librement monodromique $\nablatilde_v\in\HkGmonodrom$ sous Koszul, voir~\cite{bezrukavnikovKoszulDualityKacMoody2013}.

L'action de $\nablatilde_v$ sur $C_z$ donne tout simplement la convolution par rapport à l'Iwahori  $C_z*\nabla_v$, laquelle est filtrée par des objets simples dans des cellules à la Lusztig plus petites ou égales à celle de~$z$. Compte tenu de notre description de $\Ask$ en termes des objets simples, pour que le Hom ci-dessus soit non-nul, il faut que la cellule, ou plutôt l'orbite nilpotente~$O'$, associée à~$z$ soit plus grande ou égale à l'orbite initiale~$O$ correspondant à~$\sigma$; il faut aussi que~$z$ soit sphérique à gauche, c'est-à-dire la multiplication à gauche par tout élément du groupe de Weyl de~$G$ l'allonge. On va choisir le niveau $\kappa$ vérifiant une condition assez restrictive en plus.
\begin{cond} \label{sing}
    Dans la $O$-cellule, il existe des éléments qui sont sphériques à gauche et en même temps peuvent terminer par toutes les singularités de $\Sing(\kp)$, mais pas dans une cellule plus grande.
\end{cond}
Pour ce genre de couple $(O,\kappa)$, $\Vsk$ sera la section globale $\kp$-tordue de $\ICOk*\Deltavi$ dont la détermination est largement réduite aux calculs de $C_z*\nabla_v$ pour les~$z$ de longueur maximale dans $z\Sing(\kp)$, sphériques à gauche, et se trouvant dans la $O$-cellule. Lorsque~$\kappa$ est proche de zéro, $v$ sera relativement court, ces convolutions se calculeraient à la main. On va procéder au cas par cas, selon le type de~$G$, l'orbite~$O$, et le niveau~$\kappa$. On a déjà traité le cas de l'orbite régulière d'une manière plus directe. La cellule associée est constituée simplement des automorphismes~$\gamma$ définis plus haut. Pour les autres orbites, la condition restrictive ci-dessus requiert un $\kappa<0$.

L'orbite qui lui succède naturellement est la sous-régulière. La cellule correspondante est constituée des éléments dont l'expression réduite est unique à déplacement du facteur $\gamma$ près. Ceux qui sont déjà dans la cellule régulière sont exclus. Pour qu'un tel élément puisse terminer par toutes les singularités de~$\kp$, il faut qu'exactement une réflexion simple, notée~$t$, soit singulière pour~$\kp$. Ceci contraint le niveau à $-1,\ldots,-d$ où~$d$ est le coefficient maximal qui apparaît lorsqu'on exprime la coracine duale à la plus haute racine de~$G$ comme combinaison linéaire des coracines simples. Pour un tel niveau, l'élément~$v$ est déterminé par le fait que $tv$ est l'élément le plus court d'expression réduite unique commençant par~$t$ et terminant par~$s_0$. On peut alors définir une bijection entre les éléments~$w$ sphériques de deux côtés dans la cellule sous-régulière et les éléments~$z$ dans la même cellule qui sont sphériques à gauche mais terminent par~$t$: on associe $z=v^{-1}t$ à $w=s_0$, et $z=wv^{-1}$ aux autres~$w$. Si~$z$ et~$w$ se correspondent, $C_z*\nabla_v$ sera isomorphe à $C_w$ modulo les objets simples négligeables qui soit ne se trouvent plus dans de la cellule sous-régulière, soit ne sont pas sphérique à droite. Ceci se montre aisément à l'aide de la formule récursive classique dans~\cite{kazhdanRepresentationsCoxeterGroups1979}:\[C_r*C_s\ito C_{rs}\oplus\bigoplus_{r's<r'<r}C_{r'}^{\oplus\mu(r',r)}\]
pour toute réflexion simple $s$ et tout $r$ ne pouvant pas terminer par~$s$, où $\mu$ est le coefficient de Lusztig dont le sens est le suivant. Considérons l'extension intermédiaire du faisceau constant sur la strate $r$ vers les strates plus petites. D'après la formule l'exprimant comme une itération des images directes tronquées, on sait que les décalages cohomologiques qui en résultent sur sa tige sur un point d'une strate $r'<r$ sont toutes majorées par la codimension de la strate moins un. Le coefficient $\mu(r',r)$ compte la multiplicité de cette décalage maximale. Il pourrait être nul. 

On en déduit que $\Hom(P_z,\IC_w*\Deltavi)\ito\CC$ si $z,w$ se correspondent et est nul sinon. Il s'en suit que $\IC_w*\Deltavi$, et donc aussi $\Ask*\Deltavi$, est isomorphe à $\IC_z$ modulo les objets simples négligeables qui sont  soit constants sur la $\chi$-composante, soit ne terminent pas par~$t$, où $z$ est associé à~$w$. Ceci implique que $\Vsk$ est isomorphe à $L(z\cdot \kp)$. Rappelons que pour l'orbite sous-régulière, la représentation triviale de $\Zs$ correspond à $w=s_0$. Il résulte que l'algèbre vertex $\Vsk[1]$ est~$\Lk$, le quotient simple du vacuum. Son  module factorisable $\Vsk$ pour $\chi$ non-trivial est $L(w\cdot \kappa\Lhat_0)$. On voit que pour $\kappa=-1$, on retrouve la formule déjà obtenue.

Pour la convenance du lecteur et pour faciliter la comparaison avec la littérature, nous produisons la table~\ref{sous} rendant ces résultats explicites, en l'accompagnant d'une brève explication des notations. 

\begin{table}[htbp]
\footnotesize
  \caption{Cas sous-réguliers}\label{sous}
\[\begingroup \renewcommand{\arraystretch}{1.5}
\setlength{\arraycolsep}{0pt}
\begin{array}{lcccll}
G & \Zs & \chi & w & \multicolumn{2}{c}{\Vsk} \\
\hline
\Al,\ \ell\geqslant 2&\Gm & n\text{ si }i=1; -n\text{ sinon} & (s_0\gamma_i)^n s_0 & \multicolumn{2}{l}{\LkL[-1]{n\Lambda_j}\text{ où }\{i,j\}=\{1,\ell\}}\\
\hline
\rmA_1 &\SL_2 & n & (s_0\gamma_1)^n s_0 & \LkL[-1]{n\Lambda_1}&\\
\hline
&&&s_0  &  \Lk[-1]&\Lk[-2]\\
 &&   &s_0s_2s_1\gamma_1  & \LkL[-1]{\Lambda_3}&\Lk[-2]\\
 \Bl & (\ZZ/2)^2 &&& (\LkL[-1]{2\Lambda_3}\text{ si }\ell=3)&\\
 \ell\geqslant 3 && & s_0s_2\cdot\cdot s_{\ell}\cdot\cdot  s_2s_0 & \LkL[-1]{(2\ell-3)\Lambda_1}& \LkL[-2]{(2\ell-5)\Lambda_1}\\
  &&& s_0s_2\cdot\cdot   s_{\ell}\cdot\cdot  s_1\gamma_1 & \LkL[-1]{(2\ell-4)\Lambda_1+\Lambda_2} & \LkL[-2]{(2\ell-5)\Lambda_1}\\
  \hline
\multirow{2}{*}{$\Cl$} & \multirow{2}{*}{$\ZZ/2\ltimes\Gm$}& n & (s_0s_1\cdot\cdot  s_{\ell-1}\gamma_{\ell})^n s_0 &\LkL[-1]{n\ell\Lambda_1}\\
 &&   \text{sign.} & s_0s_1s_0 & \LkL[-1]{\Lambda_2}&\\
\hline
\multirow{5}{*}{$\Dl$}&  &&s_0 & \Lk[-1]&\Lk[-2]\\
   &\ZZ/4 & \text{sign.}&  s_0s_2s_1\gamma_1 & \LkL[-1]{\Lambda_3}&\Lk[-2]\\
   &\text{si }2\nmid\ell;&&&(\LkL[-1]{\Lambda_3+\Lambda_4}\text{ si }\ell=4) &\\
 &(\ZZ/2)^2&&    s_0s_2\cdot\cdot s_{\ell-2}s_i\gamma_i& \LkL[-1]{(\ell-3)\Lambda_1+\Lambda_j}&\LkL[-2]{(\ell-4)\Lambda_1}\\
&\text{sinon}&&&\text{où }\{i,j\}=\{\ell-1,\ell\}&\\
 \hline
 \end{array}\endgroup
 \] \vspace{-\topsep}
 \[\begingroup \renewcommand{\arraystretch}{1.5}
\setlength{\arraycolsep}{7pt}
\begin{array}{lcccllllll}
 &&&  s_0 & \Lk[-1] & \Lk[-2]& \Lk[-3]&&&\\
\rmE_6 & \ZZ/3 &  & s_0s_6s_3\gamma_i s_6s_0 & \LkL[-1]{2\Lambda_j}&\LkL[-2]{\Lambda_j}& \Lk[-3]&&&\\
   && &  &\multicolumn{6}{l}{\text{où }(i,j)=(1,4) \text{ ou }(5,2)}\\
   \hline
    \multirow{2}{*}{$\rmE_7$}& \multirow{2}{*}{$\ZZ/2$} && s_0 & \Lk[-1] &\Lk[-2] & \Lk[-3] &\Lk[-4]&&\\
   &&& s_0s_1\cdot\cdot  s_6\gamma_6 & \LkL[-1]{3\Lambda_7}&\LkL[-2]{2\Lambda_7}& \LkL[-3]{\Lambda_7} &\Lk[-4] &&\\
   \hline
   \rmE_8 &1& & s_0 & \Lk[-1]&\Lk[-2]& \Lk[-3]&\Lk[-4]&\Lk[-5]&\Lk[-6]\\
   \hline
  \multirow{2}{*}{$\rmF_4$}&\multirow{2}{*}{$\ZZ/2$}&& s_0 & \Lk[-1] &\Lk[-2]&\Lk[-3]&&&\\
  && & s_0s_1s_2s_3s_2s_1s_0 & \LkL[-1]{5\Lambda_4}&\LkL[-2]{3\Lambda_4}&\LkL[-3]{\Lambda_4}&&&\\
  \hline
 & & & s_0 & \Lk[-1]&\Lk[-2]&&&&\\
   \rmG_2 & S_3 &\text{sign.} &s_0s_1s_2s_1s_0 & \LkL[-1]{4\Lambda_2}&\LkL[-2]{\Lambda_2} &&&&\\
    &&2\text{-dim.} & s_0s_1s_2s_1s_2s_1s_0 & \LkL[-1]{3\Lambda_1}&\LkL[-2]{\Lambda_1}&&&&\\
    \hline
\end{array}\endgroup
\]
\end{table}

Dans ce tableau,~$n$ parcourt les entiers positifs. La représentation irréductible~$\chi$ indexée par $n$ est de plus haut poids~$n$. Le mot \emph{sign.} est une abréviation pour \enquote{signature}, et dénote un caractère d'ordre~2. Aux cas où des $\chi$ sont en position symétrique, on se dispense de spécifier leur correspondance exacte avec les~$w$. Notre numérotation des sommets des diagrammes de Dynkin de type fini suit la table~\emph{Fin} en page~53 du  livre~\cite{kacInfiniteDimensionalLieAlgebras1990}. On y ajoute la sommet numérotée~0 pour former le digramme de Dynkin affine.

On peut aller au-delà du cas sous-régulier et traiter des orbites~$O$ encore plus petites. Rappelons que notre condition restrictive sur le couple $(O,\kappa)$ requiert qu'on aie un contrôle sur toutes les cellules plus grandes ou égales à celle de~$O$.
Il est plus aisé de vérifier cette condition si l'on connaît explicitement tous les éléments dans ces cellules. Une notion introduite par Stembridge simplifiera cette tâche. Il appelle un élément \emph{pleinement commutatif} lorsque aucune de ses expressions réduites ne contient, comme sous-expression convexe, 
l'élément le plus longue dans un sous-groupe de Coxeter non-commutatif engendré par deux réflexions simples. De tels éléments sont relativement faciles à reconnaître et énumérer puisqu'ils admettent une interprétation visuelle, ou graphique, comme un \emph{empilement}, selon la terminologie de Viennot formalisant l'image intuitive de celui formé par des dominos sur un échiquier, voir~\cite{biagioliFullyCommutativeElements2015} pour plus de détails. Nous allons considérer des cas favorables où la condition supplémentaire suivante est satisfaite.
 \begin{cond} \label{comm}
     Toutes les cellules plus grandes ou égales à~$O$ sont constituées entièrement des éléments pleinement commutatifs.
 \end{cond}
 On voit directement que les orbites qu'on a déjà traité, à savoir la régulière et la sous-régulière, en sont des exemples. Leurs codimensions dans le cône nilpotent sont zéro et deux respectivement. Ce qui est moins connu est qu'en fait une orbite de codimension quatre satisfait aussi à cette  condition. Les inégalités~\cite[\Thmdot{Prop}{2.4}]{lusztigCellsAffineWeyl1985} et \cite[\Thmdot{Prop}{1.2}]{lusztigCellsAffineWeyl1987} permettent de le voir. On appelle ces cas sous-sous-réguliers. D'autre part, la situation pour~$G$ de type~$\Al$ ou~$\Cl$ est particulièrement simple --- Jianyi Shi a montré dans~\cite{shiFullyCommutativeElements2003,shiFullyCommutativeElements2005a} que dans ces cas, toute cellule qui contient un élément pleinement commutatif satisfait à la condition~\ref{comm}. Pour ces deux genres de situations, nous avons trouvé tous les niveaux~$\kappa$ satisfaisant à la condition~\ref{sing}. Un cas où il n'existe aucun tel niveau sera omis. La table~\ref{commAC} résume nos résultats. 
 \begin{table}[htbp] \footnotesize \caption{Cas sous-sous-réguliers, et plus} \label{commAC}
\[\begingroup\renewcommand{\arraystretch}{1.1}
\setlength{\arraycolsep}{4pt}
\begin{array}{lcccc}
G & \Zs & \chi & w & \Vsk\\
\hline
\multicolumn{3}{c}{\ker\bigl((\ell+1+\kappa,-\kappa):\Gm^2\bigr)\; (-n,0)}&\wthree{-\kappa,n-\kappa,0}&\LkL{n\Lambda_{-\kappa}}\\ 
 &\text{ si }\kappa\in [-\frac{\ell}{2},-2];&(b+n',b) & \wthree{n'-\kappa,-\kappa-b,b}&\LkL{(\ell+1+2\kappa-n')\Lambda_b+n'\Lambda_{\ell+1+\kappa+b}}\\
&\ker(\det^{\frac{\ell+1}{2}}:\GL_2)& (b-n,b)&\wthree{-\kappa,n-\kappa-b,b}&\LkL{(\ell+1+2\kappa)\Lambda_b+n\Lambda_{-\kappa}}\\
\Al&\text{ si }\kappa=-\frac{\ell+1}{2}& (n,0)&\wthree{\ell+1+\kappa,0,n-\kappa}&\LkL{n\Lambda_{\ell+1+\kappa}}\\
\multicolumn{5}{r}{(b+\ell+1+2\kappa+n,b)\;\wthree{\ell+1+\kappa,-\kappa-b,b+n}\;\LkL{n\Lambda_{\ell+1+\kappa}+(\ell+1+2\kappa)\Lambda_{\ell+1+\kappa+b}}}\\
\multicolumn{5}{l}{\text{Ici }\wthree{m,u,v}:=\bigl(01\cdots (m-1)\gamma_{\ell}\bigr)^u\bigl(0\ell\cdots(m+1)\gamma_1\bigr)^v\gamma_m, \text{ et }b\in[1,-\kappa-1], n'\in [1,\ell+2\kappa].}\\\multicolumn{5}{l}{\text{Omets les~}\chi\text{ non-dominants lorsque~}\kappa=-\frac{\ell+1}{2}.}\\
\hline\end{array}\endgroup
 \] \vspace{-\topsep}
 \[\begingroup \renewcommand{\arraystretch}{1.1}\setlength{\arraycolsep}{-10pt}
 \begin{array}{lcccc}
&&& 021320  &  \Lk[-3]\qquad\qquad\qquad\Lk[-4]\\
 \Bl & (\ZZ/2)^2 && 0213\cdots \ell\cdots 320  & \LkL[-3]{(2\ell-6)\Lambda_1} \; \LkL[-4]{(2\ell-9)\Lambda_1}\\
 \ell\geqslant 5 && & 0213240321\gamma_1 & \LkL[-3]{\Lambda_5} \ (\LkL[-3]{2\Lambda_5}\text{ si }\ell=5)\quad \Lk[-4]\\
 \multicolumn{5}{r}{021324\cdots \ell\cdots 40321\gamma_1\;\LkL[-3]{(2\ell-8)\Lambda_1+\Lambda_4}\;\LkL[-4]{(2\ell-9)\Lambda_1}}\\
  \hline
\multirow{2}{*}{$\rmB_4$}&\multirow{2}{*}{$\ZZ/2\ltimes\Gm$}&n&(0234\gamma_1)^n 021320&\LkL[-3]{2n\Lambda_4}\\&&\text{sign.}&02134320 & \LkL[-3]{\Lambda_1}\\\hline
 &\text{d'ordre~8;}  
  && \prod_{i=0}^{-2\kappa-2}(i\cdots 0)  &  \Lk \\
 &\text{diédral si}
&&\prod_{i=0}^{-2\kappa-1}(i\cdots 0) &\LkL{\Lambda_{-2\kappa}}\\
&2\mid\ell\text{ ou }2\nmid\kappa,& &
01\cdot\cdot\ell\cdot\cdot (-2\kappa-2)\prod_{i=0}^{-2\kappa-3}(i\cdot\cdot 0)& \LkL{(2\ell+3+3\kappa)\Lambda_1} \\
\Cl &\text{ hamiltonien}&\multicolumn{3}{r}{01\cdot\cdot\ell\cdot\cdot (-2\kappa-1)\prod_{i=0}^{-2\kappa-2}(i\cdot\cdot 0)\;
\LkL{(2\ell+3+4\kappa)\Lambda_1+\Lambda_{-2\kappa-1}}}\\
\multicolumn{2}{l}{\kappa\in [-\frac{\ell}{2},-2]\quad\text{sinon}}&2\text{-dim.}&(01\cdots\ell)^{-\kappa}\gamma_{\ell} & \LkL{(\ell+2+2\kappa)\Lambda_{-\kappa}} \\
\hline
\multicolumn{2}{l}{\Cl,\ 2\nmid\ell\quad\ZZ/2\ltimes\Gm\text{ si}}&2n&(01\cdots\ell)^n\prod_{i=0}^{\ell-1}(i\cdots 0)&\LkL[-\frac{\ell+1}{2}]{\Lambda_r+(n+[\frac{n}{\ell}])\Lambda_{\ell}}
\\
\multicolumn{2}{l}{\kappa=-\frac{\ell+1}{2}\quad 4\mid (\ell+1);}&2n-\ell &\multicolumn{2}{r}{(01\cdots\ell)^{n}\gamma_{\ell},\text{ ici }n\geqslant\frac{\ell+1}{2}\;\LkL[-\frac{\ell+1}{2}]{\Lambda_r+(n+[\frac{n}{\ell}]-\frac{\ell+1}{2})\Lambda_{\ell}}}\\\multicolumn{2}{r}{\Normalizer{\SL_2}{\Gm}\text{ sinon}} 
&\text{sign.}&01\cdots\ell (\ell-1)\prod_{i=0}^{\ell-2}(i\cdots 0)&\LkL[-\frac{\ell+1}{2}]{2\Lambda_1}\\\multicolumn{5}{l}{\text{Ici }N\text{ signifie le normalisateur; }r:=n-\ell[\frac{n}{\ell}];\text{ omets }\Lambda_r\text{ si }r=0.}\\\hline
\Cl,\ 2\mid\ell&\multirow{2}{*}{$\ZZ/2\times\SL_2$}&n&
(0\cdots\ell)^n\prod_{i=0}^{\ell}(i\cdots 0)&\LkL[-\frac{\ell}{2}-1]{n\Lambda_{\ell}}\\\kappa=-\frac{\ell}{2}-1&&\text{sign.}\boxtimes n&
(0\cdots\ell)^{n+\frac{\ell}{2}+1}\gamma_{\ell} &\LkL[-\frac{\ell}{2}-1]{n\Lambda_{\ell}}\\\hline
&\ZZ/4 && 021320 & \Lk[-3]\quad\Lk[-4]\\
  \Dl &\text{si }2\nmid\ell;& \text{sign.}&  0213240321\gamma_1 & \LkL[-3]{\Lambda_5}\quad\Lk[-4]\\\ell\geqslant 6
   &(\ZZ/2)^2&\multicolumn{3}{r}{02\cdot\cdot (\ell-2)j12\cdot\cdot(\ell-2)i\gamma_i\;\LkL[-3]{(\ell-5)\Lambda_2+\Lambda_i}\;\LkL[-4]{(\ell-6)\Lambda_2}}\\
 &\text{sinon}& &\text{où }\{i,j\}=\{\ell-1,\ell\}&\\
 \hline
 \rmD_5 &\Gm &n\text{ si }i=4; -n\text{ sinon}& 02(13\gamma_i)^n 1320&\LkL[-3]{n\Lambda_i}\text{ où }i=4,5\\
 \hline
\rmE_6 & \Gm & n\text{ si }i=1; -n\text{ sinon} & 063243(\gamma_i i243)^n60 & \LkL[-4]{2n\Lambda_i}\text{ où }i=1,5 \\
    \hline
    \multirow{2}{*}{$\rmE_7$}& \multirow{2}{*}{$\ZZ/2$} && 0123473210 & \Lk[-5]\quad\Lk[-6] \\
   &&& 012347352413273456\gamma_6 & \LkL[-5]{\Lambda_1+3\Lambda_6}\quad\LkL[-6]{2\Lambda_6}\\
   \hline
   \rmE_8 &1& & 01234568543210  & \Lk[-7]\quad\Lk[-8]\quad\Lk[-9]\quad\Lk[-10]\\
   \hline
  \multirow{2}{*}{$\rmF_4$}&\multirow{2}{*}{$\ZZ/2$}&& 0123243210 & \Lk[-4]\\
  && & 012324130213243210 & \LkL[-4]{4\Lambda_4}\\
  \hline
\end{array}\endgroup\]
\end{table}
Dans ce tableau, $n$ parcourt les entiers positifs sauf si l'on précise autrement. Le paramètre $\sigma:\SL_2\to\LG$ correspond aux partitions suivantes selon le type de~$G$ lorsque celui-ci est classique:
\[\begin{array}{c|cl}
      G& \multicolumn{2}{l}{\text{partition correspondant à }\sigma}\\\hline
     \Al& \ell+1+\kappa,-\kappa &\\\hline
     \Bl& 2\ell-4,4&\\\hline
     \Cl&\ell,\ell,1 &\text{ si }\ell=-2\kappa-2\\
     &2\ell+2\kappa+1,-2\kappa-1,1&\text{ sinon}\\ \hline
     \Dl& 2\ell-5,5&
\end{array}\]
Pour que les éléments, souvent assez longues, du groupe de Weyl étendu soient plus lisibles, on ne note les réflexions simples que par leurs indices. L'expression $i\dotsm j$ désigne le produit des réflexions simples d'indices consécutives de~$i$ à~$j$. Les indices sont croissants si $i\leqslant j$ et décroissants sinon.
Les conditions~\ref{sing} et~\ref{comm} garantissent l'existence d'un élément pleinement commutatif pouvant terminer par toutes les singularités de~$\kp$. Il s'en suit que les réflexions simples singulières pour~$\kp$ commutent entre eux. Ce fait nous permet de définir une bijection entre les éléments~$w$ sphériques de deux côtés dans la $O$-cellule et les éléments~$z$ dans la même cellule qui sont sphériques à gauche mais peuvent terminer par toutes les singularités de~$\kp$. On associe $z=v^{-1}\prod_{s\in\Sing(\kp)} s$ à~$w$ lorsque celui-ci est l'involution de Duflo, autrement dit lorsqu'il correspond au~$\chi$ trivial, et associe $z=wv^{-1}$ aux autres~$w$, à une exception près: on associe $z=0213\dotsm\ell\dotsm 42$ à $w= 0213\dotsm \ell\dotsm 320$ lorsque le couple $(G,\kappa)$ est $(\Bl,-4)$ où $\ell\geqslant 5$. Si~$z$ et~$w$ se correspondent, $C_z*\nabla_v$ sera isomorphe à $C_w$ modulo les objets simples négligeables. Ceci implique que $\Vsk$ est isomorphe à $L(z\cdot \kp)$. Il résulte que l'algèbre vertex $\Vsk[1]$ est~$\Lk$, le quotient simple du vacuum. Son  module factorisable $\Vsk$ pour $\chi$ non-trivial est $L(w\cdot \kappa\Lhat_0)$ à l'exception ci-dessus près. En fait, si l'on insiste à regarder $L(w\cdot \kappa\Lhat_0)$ pour ce cas exceptionnel, on se rendra immédiatement compte qu'il n'est même pas $G[[t]]$-intégrable. Enfin, le calcul explicite de la bijection de Lusztig-Vogan pour type~$\Al$ est le sujet de la thèse d'Achar dont une partie est publiée dans~\cite{acharEquivariantKtheoryNilpotent2004}, à laquelle nous nous référons. Tout comme ce qu'on a fait pour les types~$\Al,\Cl$, on espère qu'on pourrait traiter les niveaux plus généraux pour les types~$\Bl,\Dl$. Notre méthode permet déjà un calcul des modules éventuels dans ces cas et de leurs lois de fusions. Pour conclure, il faudrait encore une meilleure connaissance des cellules sous-jacentes.

Nous n'avons imposé la condition sur la pleine commutativité que pour avoir une connaissance explicite des éléments dans les cellules concernées. Parfois, la même connaissance parfaite pourrait être atteinte sans cette condition. Ce scénario alternatif est illustré par la table suivante incluant tous les couples $(O,\kappa)$ vérifiant la condition~\ref{sing} pour~$G$ de rang~2, les couples déjà traités plus haut étant omis.
\[\begin{array}{lcccc}
G\qquad\quad O&\Zs & \chi & w & \Vsk\\
\rmA_2\qquad\text{ zéro}&\SL_3&m\Lambda_1+n\Lambda_2&(01\gamma_2)^{n+1}(02\gamma_1)^{m+1}0\quad&\LkL[-2]{m\Lambda_1+n\Lambda_2} \\
\rmG_2\quad\text{minimale}&\SL_2&n&01212010(210)^n21210\quad&\LkL[-3]{n\Lambda_1}\text{ si }2\mid n\\&&&\multicolumn{2}{r}{\LkL[-3]{n\Lambda_1+\Lambda_2}\text{ sinon}}
\end{array}\]
Les éléments dans ces cellules se voient dans les dessins colorés dans~\cite{lusztigCellsAffineWeyl1985}. Les bijections de Lustig-Vogan dans ces cas sont décrites dans~\cite{xiRepresentationsAffineHecke1994}. Les $n,m$ dans ce tableau parcourent les entiers positifs. La bijection entre les~$z$ et les~$w$ est simplement donnée par~$z=wv^{-1}$.

Laissons-nous terminer par un survol des littératures étudiant soit les variétés associées des algèbres vertex affines simples à un niveau entier non-admissible, soit les modules sur ces algèbres. La table~\ref{tab:Comparaison} rassemble celles connues à l'auteur. Le lecteur est invité à les consulter pour plus de détails. 

\begin{table}[htb]\caption{Littératures liées}
\label{tab:Comparaison}
\setlength{\tabcolsep}{3pt}
\renewcommand{\arraystretch}{1.2}
\centering
\begin{tabular}{lcrr}
$G$ & $\kappa$ & modules &  variété associée\\\hline
$\Al,\ell\geqslant 2$    & -1 &\cite[6.1]{adamovicFusionRulesComplete2014} & cas $\ell=2$, nappe \cite[1.2]{jiangAssociatedVarietiesSimple2024}\\
&&&cas $\ell\geqslant 3$, nappe \cite[1.1]{arakawaSheetsAssociatedVarieties2017}\\
&$-\frac{\ell+1}{2}$ &&nappe \cite[1.1]{arakawaSheetsAssociatedVarieties2017}\\\hline
$\Bl$& -2 & \cite[7.3]{adamovicApplicationCollapsingLevels2020}*  &racine courte \cite[7.1]{arakawaIrreducibilityAssociatedVarieties2018} \\&$1-\ell$ & \multicolumn{2}{l}{\cite[4]{adamovicGeneralResultsConformal2013}*; cas $\ell=3$ \cite[C]{arakawaSimpleOrbifoldAffine2024}}\\\hline
$\Cl,\ell\geqslant 3$& -1 &\cite[6.1]{adamovicFusionRulesComplete2014}* 
&\\ 
$\rmC_2$& -2 & \cite[7.4]{adamovicApplicationCollapsingLevels2020}  &\\\hline
$\Dl$& -1&  &minimal \cite[1.1]{arakawaJosephIdealsLisse2018}\\\ 
& -2 & \cite[6.4]{adamovicApplicationCollapsingLevels2020}* &minimal \cite[1.1]{arakawaJosephIdealsLisse2018}\\
&&\multicolumn{2}{l}{cas $\ell=4$ \cite[4.3]{perseNoteRepresentationsAffine2013}, \cite[3.1]{arakawaJosephIdealsLisse2018}}\\
&$2-\ell$&cas $2\nmid\ell$, \cite[7.2]{adamovicFusionRulesComplete2014}*&nappe \cite[1.2]{arakawaSheetsAssociatedVarieties2017}\\
&&cas $2\mid\ell$, \cite[1.1]{adamovicApplicationCollapsingLevels2020}*&$\overline{\mathbb{O}_{(2^{\ell-2},1^4)}}$ \cite[6.1]{arakawaIrreducibilityAssociatedVarieties2018}\\\hline
$\El$ & $-m,\dotsc,-1$ &cas $\kappa=-m$ \cite[3.1]{arakawaJosephIdealsLisse2018} & minimal \cite[1.1]{arakawaJosephIdealsLisse2018}\\&&où $m=3,4,6$ pour $\ell=6,7,8$&\\
$\rmE_8$ &-10&\cite[1.1]{adamovicApplicationCollapsingLevels2020}&\\\hline
$\rmF_4$&-3 & \cite[4]{adamovicGeneralResultsConformal2013} &\\ \hline
$\rmG_2$ & -1 && minimal \cite[1.6]{arakawaSheetsAssociatedVarieties2017}\\ 
& -2 &\cite[5]{adamovicGeneralResultsConformal2013}, \cite[B]{arakawaSimpleOrbifoldAffine2024}& sous-régulier \cite[A]{arakawaSimpleOrbifoldAffine2024}\\\hline
\end{tabular}
\end{table}

Ces articles sont focalisés plus sur la classification et utilisent des outils très différents de la nôtre. Ils utilisent notamment les $W$-algèbres et les vecteurs singuliers. Un article marqué d'une étoile contient plus de modules, sous-entendus $G[[t]]$-intégrables, que nous. Une part de cet écart n'est essentielle, étant due à notre compréhension actuelle insuffisante des cellules des types $\Bl,\Dl$. Les \enquote{nappes} sont celles de Dixmier. Elles ne sont pas contenues dans le cône nilpotent. Pour retrouver, au sens de Spaltenstein, le dual de l'orbite nilpotente~$O$, au lieu de regarder la variété associée pour l'algèbre vertex simple, il faudrait considérer celle pour son extension~$\Vs$ --- l'algèbre de Hecke chirale nilpotente que nous avons introduite ici.

Les fusions entre ces modules étaient très peu discutées: on trouve \cite[6.2]{adamovicFusionRulesComplete2014} pour les modules sur~$\Lk[-1](\SL_n)$ et \cite[4.15]{adamovicSemisimplicityModuleCategories2024} pour ceux sur $\Lk[-2](\rmG_2)$, tous deux déjà mentionnés au début de cet article.
\printbibliography
\end{document}